\documentclass[a4paper, BCOR=0.0mm, DIV=calc, halfparskip, fleqn]{scrartcl}
\usepackage[T1]{fontenc}
\usepackage[utf8]{inputenc}
\usepackage{lmodern}
\usepackage{microtype}
\usepackage[english]{babel}
\usepackage{amsmath, amssymb, amsfonts}
\usepackage{nicefrac}
\newcommand{\diff}[1]{\mathrm{#1}}
\newcommand{\gam}[1]{\Gamma( #1 )}
\newcommand{\pham}[1]{\left( #1 \right)}
\newcommand{\nfrac}[2]{\nicefrac{#1}{#2}}

\KOMAoptions{twoside=false, twocolumn=false, headinclude=false, footinclude=false, mpinclude=false, pagesize=auto}
\recalctypearea

\begin{document}
\title{\boldmath On proving some of Ramanujan's formulas for $\frac{1}{\pi}$ with an elementary method\unboldmath}
\author{Alexander Aycock}
\date{}
\maketitle
\begin{abstract}
In this paper we want to prove some formulas listed by S. Ramanujan in his paper "Modular equations and approximations to $\pi$" \cite{24} with an elementary method.
\end{abstract}
\section*{Introduction}
In 1914 S. Ramanujan publishes his nowaday famous memoir "Modular equations and approximations to $\pi$" \cite{24}, in which he gives, along with sixteen other formulas of the same type, this particular formula:
\[
\frac{1}{2\pi\sqrt{2}}=\frac{1103}{99^2}+\frac{27493}{99^6}\frac{1}{2}\frac{1\cdot3}{4^2}+\frac{53883}{99^{10}}\frac{1\cdot3}{2\cdot4}\frac{1\cdot3\cdot5\cdot7}{4^2\cdot8^2}+\cdots
\]
Ramanujan gives the formula as seen above, and in his treatise it is formulas (44). This formula is not only memorable because of its nature, but also, because it converges rapidly and is appropiate, to calculate decimal digits of $\pi$. It gives, as also Ramanujan notes, by taking the first term alone
\[
\frac{1103}{99^2}=0,1125395678...
\]
and in comparision we have
\[
\frac{1}{2\pi\sqrt{2}}=0.11253953951...
\]
so that the values differ at the ninth decimal.
How exactly Ramanujan arrived at his formulas, is, even by studying his paper intensively, not easy to figure out and so it comes, that he did not prove his formulas with all rigour.\\

They indeed remained unproved until the paper of J. Borwein und P. Borwein "More Ramanujan-type series for $\frac{1}{\pi}$" \cite{5}, which was published in 1987. The proofs are given by using the properties the $j-$function, Klein's absolute invariant, about which we will not talk about here any further, because our method of proof will be a different one, which we will explain later in all detail.\\

But at first we want to note in this intrductory part, that along with the idea involving $j-$function, as also D. Chudnovsky und G. Chudnovsky in their paper "Approximations and complex multiplication according to Ramanujan" \cite{9} use it, to prove more formulas of this type, including this one
\[
	\frac{1}{\pi} = 12\sum_{n=0}^{\infty} \frac{(-1)^k (6k)! (15591409+545140134k)}{(3k)!(k)!^3 640320^{3k+\frac{3}{2}}}	
\]
which because of its rapid convergence is also quite appropriate for calculating decimal places of $\pi$, there are other methods of proof for this formulas. Here we want to mention the WZ-Algorithm, with which D. Zeilberger in the paper "A WZ proof of Ramanujan's formula for $\pi$" \cite{30} gives a proof for a dertain formula. Also J. Guillera often used this algorithm for the same purpose, and proved this analog series $\frac{1}{\pi^2}$ for example
\[
	\sum_{n=0}^{\infty} (13 + 180n + 820n^2) \frac{(\frac{1}{2})_n^5}{(1)_n^5}\frac{(-1)^n}{2^{10}} = \frac{128}{\pi^2}
\]
where - as constantly in the following -
\[
(a)_n=\frac{\Gamma{(a+n)}}{\Gamma{(a)}}=\prod_{k=1}^n(a+k)
\]
with the additional condition
\[
(a)_0=1
\]
is the Pochhammer symbol, $\Gamma{(x)}$ is the Gamma function, defined as
\[
\Gamma{(x)}=\int_0^{\infty}e^{-t}t^{x-1}\diff{d}t
\]
But we will also not pursue this way any further, but go another way, which is likely to be seen as more elementary as the two already mentionend - also for that reason, because it requires less work in advance and nevertheless leads to some (although not to all) formulas for $\frac{1}{\pi}$, that Ramanujan recorded in his paper \cite{24}. And to the explaination and the application of this method, to get to the Ramanujan formulas for $\frac{1}{\pi}$, this paper will be devoted. Although if one already finds examples of this methods scattered around everywhere - W. Zudilin in his paper "Lost in Translation" \cite{31} gives some references and explains the method on an example, to which we will retutn later -, you find no memoir, which solely uses this method. Instead you have to use a already known formula, to get a new one from it. This already sums up the nature of the translation method.\\

The translation method is a method, to prove from an equation $A=B$ by means of transformations, substitutions etc. a "new" one $C=D$.\\

The word "new" was highlighted on purpose, because the equations $A=B$ and $C=D$ are equivalent and therefore $C=D$ is actually not new. And vice versa you could prove $A=B$ with the same transformation, substitution etc. just used in the other direction, from $C=D$.\\

This at the same time reveals a weakness of this method, even if the translation can be done by elementary means, one can nevertheless only show the equivalence of both equations, but not their validity. One of the two equations must be treated as known, which is likely the reason, that there is no treatise, only relying on this method. You need a formula, you can start from.\\  

And for our undertaking, the proof of some formulas $\frac{1}{\pi}$, that Ramanujan gave in \cite{24}, it is even possible to find such a starting formula, as we will see below. These things, of course, immedeatly give the structure of this paper.\\

\section*{Contens of this paper}

We mainly want to deal with the formulas for $\frac{1}{\pi}$, that Ramanujan gave, and want to prove them with the translation method, which we will explain on the concret example of these formulas.\\

At first we want to give a list with formulas, that will be in the focus of our interest. Ramanujan gave 17 formulas of this kind in his paper \cite{24}, but today more are known, and J. Guillera gave a list on his website "Tables of Ramanujan series with rational values of $z$" \cite{18}. What rational means in this context, will become clear later.\\

Next we will give those things, you need to know in advance for the following investigation. This will be mainly about the Gau\ss ian hypergeometric series - $_2F_1(a,b,c,z)$ - and the transformation formulas and identities, it satisfies and we will make use of, to prove at least some formulas for $\frac{1}{\pi}$, that also Ramanujan gave, by their means, after we will have derived a starting formula with well-known theorems about the hypergeometric series .\\

But we will realize, that we will not be able, to prove all formulas from the mentioned list by the translation method. The reason for this will be, that we will run out of the transformations for the hypergeometric series. So we want to see at first, how far we can get, without using any other formulas from the list. Then we want to assume the formulas given by Ramanujan as known and see, how far Ramanujan's list could be extended from there - and could have been extended by Ramanujan himself. After this we want to do this with the formulas, proved by the WZ-method.\\

Because the hypergeometric series is a power series in its last variable, considerations about its convergence are neccessary, of course, but nevertheless will the translation method also lead to divergent series of the corresponding shape - these will also be found on the list -, which will have a meaning, that we will explain then.\\

After having treated the Ramanujan-formulas as good as possible with the translation method, we will consider formulas of the same shape, that arise from limit formulas, formulas of which kind seem to be hardly considered up to now. At this point we will treat all formulas from our list as known, to get more possibilties.\\

At last we want to apply the translation method to some chosen examples, given by Z.W. Sun in the paper "List of conjectural series for powers of $\pi$ and other constants" \cite{29} and want to prove the corresponing conjectures.\\

Then, being content, to have explained the translation method on these examples, we want to conclude and give a short review and overview.\\

So it follows the list of formulas, that will be of interest for us.\\

\section*{List of Ramanujan series for rational values of $z$}

All series have this form
\[
\sum_{n=0}^{\infty}\frac{(\frac{1}{2})_n(s)_n(1-s)_n}{(1)_n^3}(a+bn)z^n=\frac{c}{\pi}
\]
where $s=\frac{1}{2},\frac{1}{3},\frac{1}{4}$ or $\frac{1}{6}$ , $z$ is a rational value, $a$ and $b$ are positive whole numbers, $c^2$ is a rational number.\\
Of course all these restrictions are not neccessary, but these cases are the most convinient ones (easy numerical evaluation, calculation etc.) and they will suffice, to illustrate the translation method.\\

In the following list the "R" indicates, that the corresponding formula was already given by Ramanujan himself. "WZ" indicates, that the formula was proved by the WZ-method. If there is no letter, this means, that the formula was proved by other means, in most cases by modular forms. We divide the series according to the values of $s$ - following Guillera - into classes. \\

\section*{The list}
\subsubsection*{\boldmath $s = \frac{1}{2}$ \unboldmath}
\begin{flalign*}
	& \sum_{n=0}^{\infty} \frac{(\frac{1}{2})_n^3 }{(1)_n^3}(4n + 1)\pham{-1}^n = \frac{2}{\pi} &  WZ, (R)\\
	& \sum_{n=0}^{\infty} \frac{(\frac{1}{2})_n^3 }{(1)_n^3}(6n + 1)\pham{-\frac{1}{8}}^n = \frac{2\sqrt{2}}{\pi} & WZ & \\
	& \sum_{n=0}^{\infty} \frac{(\frac{1}{2})_n^3 }{(1)_n^3}(6n + 1)\pham{\frac{1}{4}}^n = \frac{4}{\pi} &  WZ, R &\\
	& \sum_{n=0}^{\infty} \frac{(\frac{1}{2})_n^3 }{(1)_n^3}(42n + 5)\pham{\frac{1}{64}}^n = \frac{16}{\pi} & WZ, R &\\
	&  \sum_{n=0}^{\infty} \frac{(\frac{1}{2})_n^3 }{(1)_n^3}(3n + 1)
\pham{-8}^n = \frac{1}{\pi} & WZ &\\
    &  \sum_{n=0}^{\infty} \frac{(\frac{1}{2})_n^3 }{(1)_n^3}(3n + 1)
\pham{4}^n = \frac{2i}{\pi} & WZ &\\
&  \sum_{n=0}^{\infty} \frac{(\frac{1}{2})_n^3 }{(1)_n^3}(21n + 8)
\pham{64}^n = \frac{2i}{\pi} & WZ &
\end{flalign*}
\subsubsection*{\boldmath $s = \frac{1}{4}$ \unboldmath}
\begin{flalign*}
& \sum_{n=0}^{\infty} \frac{(\frac{1}{4})_n (\frac{3}{4})_n (\frac{1}{2})_n}{(1)_n^3}(20n + 3)\pham{-\frac{1}{4}}^n = \frac{8}{\pi} & WZ, R & \\
& \sum_{n=0}^{\infty} \frac{(\frac{1}{4})_n (\frac{3}{4})_n (\frac{1}{2})_n}{(1)_n^3}(65n + 8)\pham{-\frac{16^2}{63^2}}^n = \frac{9\sqrt{7}}{\pi} &  R & \\
& \sum_{n=0}^{\infty} \frac{(\frac{1}{2})_n (\frac{3}{4})_n (\frac{1}{4})_n}{(1)_n^3}(28n + 3)\pham{-\frac{1}{48}}^n = \frac{16}{\sqrt{3}\pi} & WZ, R & \\
& \sum_{n=0}^{\infty} \frac{(\frac{1}{2})_n (\frac{3}{4})_n (\frac{1}{4})_n}{(1)_n^3}(260n + 23)\pham{-\frac{1}{18^2}}^n = \frac{72}{\pi} & R & \\
& \sum_{n=0}^{\infty} \frac{(\frac{1}{2})_n (\frac{3}{4})_n (\frac{1}{4})_n}{(1)_n^3}(644n + 41)\pham{-\frac{1}{5 \cdot 72^2}}^n = \frac{288}{\sqrt{5}\pi} & R & \\
& \sum_{n=0}^{\infty} \frac{(\frac{1}{2})_n (\frac{3}{4})_n (\frac{1}{4})_n}{(1)_n^3}(21460n + 1123)\pham{-\frac{1}{882^2}}^n = \frac{3528}{\pi} & R & \\
& \sum_{n=0}^{\infty} \frac{(\frac{1}{2})_n (\frac{3}{4})_n (\frac{1}{4})_n}{(1)_n^3}(7n + 1)\pham{\frac{32}{81}}^n = \frac{9}{2\pi}   & \\
& \sum_{n=0}^{\infty} \frac{(\frac{1}{2})_n (\frac{3}{4})_n (\frac{1}{4})_n}{(1)_n^3}(8n + 1)\pham{\frac{1}{9}}^n = \frac{2\sqrt{3}}{\pi} & WZ, R & \\
& \sum_{n=0}^{\infty} \frac{(\frac{1}{2})_n (\frac{3}{4})_n (\frac{1}{4})_n}{(1)_n^3}(40n + 3)\pham{\frac{1}{49^2}}^n = \frac{49\sqrt{3}}{9\pi} & R & \\
& \sum_{n=0}^{\infty} \frac{(\frac{1}{2})_n (\frac{3}{4})_n (\frac{1}{4})_n}{(1)_n^3}(280n + 19)\pham{\frac{1}{99^2}}^n = \frac{18\sqrt{11}}{\pi} & R & \\
& \sum_{n=0}^{\infty} \frac{(\frac{1}{2})_n (\frac{3}{4})_n (\frac{1}{4})_n}{(1)_n^3}(26390n + 1103)\pham{\frac{1}{99^4}}^n = \frac{9801\sqrt{2}}{\pi} & R & \\
& \sum_{n=0}^{\infty} \frac{(\frac{1}{2})_n (\frac{3}{4})_n (\frac{1}{4})_n}{(1)_n^3}(10n + 1)\pham{\frac{1}{81}}^n = \frac{9\sqrt{2}}{\pi} & R &\\
& \sum_{n=0}^{\infty} \frac{(\frac{1}{2})_n (\frac{3}{4})_n (\frac{1}{4})_n}{(1)_n^3}(5n + 1)\pham{-\frac{16}{9}}^n = \frac{\sqrt{3}}{\pi} & WZ &\\ 
& \sum_{n=0}^{\infty} \frac{(\frac{1}{2})_n (\frac{3}{4})_n (\frac{1}{4})_n}{(1)_n^3}(35n + 8)\pham{\frac{256}{81}}^n = \frac{18i}{\pi} & WZ &
\end{flalign*}
\subsubsection*{\boldmath $s = \frac{1}{3}$ \unboldmath}
\begin{flalign*}
& \sum_{n=0}^{\infty} \frac{(\frac{1}{2})_n (\frac{1}{3})_n (\frac{2}{3})_n}{(1)_n^3}(5n + 1)\pham{\frac{-9}{16}}^n = \frac{4}{\sqrt{3}\pi} & WZ & \\
& \sum_{n=0}^{\infty} \frac{(\frac{1}{2})_n (\frac{1}{3})_n (\frac{2}{3})_n}{(1)_n^3}(51n + 7)\pham{\frac{-1}{16}}^n = \frac{12\sqrt{3}}{\pi} & WZ & \\
& \sum_{n=0}^{\infty} \frac{(\frac{1}{2})_n (\frac{1}{3})_n (\frac{2}{3})_n}{(1)_n^3}(9n + 1)\pham{\frac{-1}{80}}^n = \frac{4\sqrt{3}}{\sqrt{5}\pi} & & \\
& \sum_{n=0}^{\infty} \frac{(\frac{1}{2})_n (\frac{1}{3})_n (\frac{2}{3})_n}{(1)_n^3}(1230n + 106)\pham{\frac{-1}{2^{10}}}^n = \frac{192\sqrt{3}}{\pi} &  & \\
& \sum_{n=0}^{\infty} \frac{(\frac{1}{2})_n (\frac{1}{3})_n (\frac{2}{3})_n}{(1)_n^3}(330n + 26)\pham{\frac{-1}{3024}}^n = \frac{216}{\sqrt{7}\pi} & & \\
& \sum_{n=0}^{\infty} \frac{(\frac{1}{2})_n (\frac{1}{3})_n (\frac{2}{3})_n}{(1)_n^3}(6n + 1)\pham{\frac{1}{2}}^n = \frac{3\sqrt{3}}{\pi} & WZ & \\
& \sum_{n=0}^{\infty} \frac{(\frac{1}{2})_n (\frac{1}{3})_n (\frac{2}{3})_n}{(1)_n^3}(15n + 2)\pham{\frac{2}{27}}^n = \frac{27}{4\pi} & R & \\
& \sum_{n=0}^{\infty} \frac{(\frac{1}{2})_n (\frac{1}{3})_n (\frac{2}{3})_n}{(1)_n^3}(33n + 4)\pham{\frac{4}{125}}^n = \frac{15\sqrt{3}}{\pi} & R & \\
& \sum_{n=0}^{\infty} \frac{(\frac{1}{2})_n (\frac{1}{3})_n (\frac{2}{3})_n}{(1)_n^3} (14151n + 827)\pham{-\frac{1}{500^2}}^n = \frac{1500\sqrt{3}}{\pi} & &\\
& \sum_{n=0}^{\infty} \frac{(\frac{1}{2})_n (\frac{1}{3})_n (\frac{2}{3})_n}{(1)_n^3} (15n + 4)\pham{-4}^n = \frac{3\sqrt{3}}{\pi} & WZ &\\
& \sum_{n=0}^{\infty} \frac{(\frac{1}{2})_n (\frac{1}{3})_n (\frac{2}{3})_n}{(1)_n^3} (10n + 3)\pham{\frac{27}{2}}^n = \frac{10i}{\pi} & WZ &\\
& \sum_{n=0}^{\infty} \frac{(\frac{1}{2})_n (\frac{1}{3})_n (\frac{2}{3})_n}{(1)_n^3} (11n + 3)\pham{\frac{27}{16}}^n = \frac{12i}{\pi} & WZ &\\
\end{flalign*}
\subsubsection*{\boldmath $s = \frac{1}{6}$ \unboldmath}
\begin{flalign*}
& \sum_{n=0}^{\infty} \frac{(\frac{1}{6})_n (\frac{1}{2})_n (\frac{5}{6})_n}{(1)_n^3}(63n + 8)\pham{\frac{-4^3}{5^3}}^n = \frac{5\sqrt{15}}{\pi} & WZ  & \\
& \sum_{n=0}^{\infty} \frac{(\frac{1}{6})_n (\frac{1}{2})_n (\frac{5}{6})_n}{(1)_n^3}(154n + 15)\pham{\frac{-3^3}{8^3}}^n = \frac{32\sqrt{2}}{\pi} & WZ & \\
& \sum_{n=0}^{\infty} \frac{(\frac{1}{6})_n (\frac{1}{2})_n (\frac{5}{6})_n}{(1)_n^3}(342n + 25)\pham{\frac{1}{8^3}}^n = \frac{32\sqrt{6}}{\pi} &  & \\
& \sum_{n=0}^{\infty} \frac{(\frac{1}{6})_n (\frac{1}{2})_n (\frac{5}{6})_n}{(1)_n^3}(506n + 31)\pham{\frac{-9}{40^3}}^n = \frac{160\sqrt{30}}{9\pi} &  & \\
& \sum_{n=0}^{\infty} \frac{(\frac{1}{6})_n (\frac{1}{2})_n (\frac{5}{6})_n}{(1)_n^3}(5418n + 263)\pham{\frac{-1}{80^3}}^n = \frac{1}{2}\cdot\frac{80^2}{\sqrt{15}\pi} &  & \\
& \sum_{n=0}^{\infty} \frac{(\frac{1}{6})_n (\frac{1}{2})_n (\frac{5}{6})_n}{(1)_n^3}(261702n + 10177)\pham{\frac{-1}{440^3}}^n = \frac{3\cdot 440^2}{\sqrt{330}\pi} &  & \\
& \sum_{n=0}^{\infty} \frac{(\frac{1}{6})_n (\frac{1}{2})_n (\frac{5}{6})_n}{(1)_n^3}(28n + 3)\pham{\frac{3}{5}}^{3n} = \frac{5\sqrt{5}}{\pi} &  \\
& \sum_{n=0}^{\infty} \frac{(\frac{1}{6})_n (\frac{1}{2})_n (\frac{5}{6})_n}{(1)_n^3}(256n + 20)\pham{\frac{2}{11}}^{3n} = \frac{11\sqrt{33}}{\pi} &  \\
& \sum_{n=0}^{\infty} \frac{(\frac{1}{6})_n (\frac{1}{2})_n (\frac{5}{6})_n}{(1)_n^3}(22n + 2)\pham{\frac{4}{5^3}}^n = \frac{5\sqrt{15}}{3\pi} & R & \\
& \sum_{n=0}^{\infty} \frac{(\frac{1}{6})_n (\frac{1}{2})_n (\frac{5}{6})_n}{(1)_n^3}(2394n + 144)\pham{\frac{4^3}{85^3}}^n = \frac{85\sqrt{85}}{\sqrt{3}\pi} & WZ, R  & \\
& \sum_{n=0}^{\infty} \frac{(\frac{1}{6})_n (\frac{1}{2})_n (\frac{5}{6})_n}{(1)_n^3}(545140134n + 13591409)\pham{-\frac{1}{53360^3}}^n= \frac{3}{2}\cdot\frac{53360^3}{\sqrt{10005}} \cdot \frac{1}{\pi}& &
\end{flalign*}
\subsubsection*{Remark}
The R in the first formula is put in brakets, because Ramanujan already knew this particular series, as it is clear from a letter from Ramanujan to G.H Hardy. See the book "Ramanujan: Letters and Commentary" by. B. Berndt and R. Rankin \cite{4}. Nevertheless this formula is not found on Ramanujan's list in \cite{24}.
\section*{Preparations - Things we need to know for the further investigation}
As already mentioned, all series on the list have this form
\[
\sum_{n=0}^{\infty}\frac{(\frac{1}{2})_n(s)_n(1-s)_n}{(1)_n^3}(a+bn)z^n=\frac{c}{\pi}
\]
with the certain conditions concerning $s, a, b, z$ and $c$.\\

Now the question remains, where this series come from. At first it becomes clear, that you only have to consider this series
\[
\sum_{n=0}^{\infty}\frac{(\frac{1}{2})_n(s)_n(1-s)_n}{(1)_n^3}z^n
\]
because this series
\[
\sum_{n=0}^{\infty}\frac{(\frac{1}{2})_n(s)_n(1-s)_n}{(1)_n^3}nz^n
\]
arises from the first, by differentiating the first with respect to $z$ and multiplying with $z$ afterwards.\\
And this series is easyly seen to be just a special case of this one
\[
G(a,b,c,d,e,z)=\sum_{n=0}^{\infty}\frac{(a)_n(b)_n(c)_n}{(d)_n(e)_n}\frac{z^n}{n!}
\]
And this series (along with other analog ones) after the first investigations by C. Gau\ss{}  "Disquisitiones generales circa seriem infinitam $1+\frac{\alpha\beta}{1\cdot\gamma}x+\frac{\alpha(\alpha+1)\beta(\beta+1)}{1\cdot2\cdot\gamma(\gamma+1)}xx+\frac{\alpha(\alpha+1)(\alpha+2)\beta(\beta+1)(\beta+2)}{1\cdot2\cdot3\cdot\gamma(\gamma+1)(\gamma+2)}x^3+\text{etc.}$" \cite{16} and L. Euler "Specimen transformationis singularis serierum" \cite{14} is called a hypergeometric series. Considering the title of Gau\ss' paper \cite{16} one sees, that the function $G(a,b,c,d,e,z)$ is a generalisation of the series investigated by Gau\ss 
\[
F(\alpha,\beta,\gamma,z)=1+\frac{\alpha\beta}{1\cdot\gamma}x+\frac{\alpha(\alpha+1)\beta(\beta+1)}{1\cdot2\cdot\gamma(\gamma+1)}xx+\frac{\alpha(\alpha+1)(\alpha+2)\beta(\beta+1)(\beta+2)}{1\cdot2\cdot3\cdot\gamma(\gamma+1)(\gamma+2)}x^3+\text{etc.}
\]
or in shorter form
\[
F(\alpha,\beta,\gamma,z)=\sum_{n=0}^{\infty}\frac{(\alpha)_n(\beta)_n}{(\gamma)_n}\frac{z^n}{n!}
\]
The series $F(\alpha,\beta,\gamma,z)$ is - because of the work mentioned by Gau\ss{} \cite{14} - often called Gau\ss ian hypergeometric series, although it was already written down by Euler in this form in the treatise cited above \cite{14} and was then investigated by him. So the name "Gau\ss ian hypergeometric series" is not - historically speaking - quite right.\\

Often the series $F(\alpha,\beta,\gamma,z)$ is also written as $_2F_1(\alpha,\beta,\gamma,z)$ and the series, we denoted by $G(a,b,c,d,e,z)$, in analogy by $_3F_2(a,b,c,d,e,z)$. The indices are self-explainatory and correspond to the number of variables in the numerator - where $z$ is not counted - and denominator of the summand. That these series can be further generalized and you can define functions as $_4F_3$, $_5F_4$, $_2F_7$ and in general $_pF_q$, where $p$ and $q$ are natural numbers, is obvious and was naturally done, of course. But for us the series $_3F_2$ and $_2F_1$ will mostly suffice, which is, why we only want to consider these ones. We will write only $F$ for $_2F_1$ and/or $_3F_2$, if this will cause no confusion\\

To use the translation method, we see, that there have to consist transformations between the the classes of the formulas for $\frac{1}{\pi}$ and therefore transformations between the corresponding hypergeometric series. In general they would have to look like this, that can also be divided into the following classes:
\begin{align*}
\text{1st Class} \qquad & F(\tfrac{1}{2}, \tfrac{1}{2}, \tfrac{1}{2}, 1, 1, A(x)) = B(x)F(\tfrac{1}{2}, \tfrac{1}{2}, \tfrac{1}{2}, 1, 1, C(x)) \\
\text{2nd Class} \qquad & F(\tfrac{1}{2}, \tfrac{1}{2}, \tfrac{1}{2}, 1, 1, A(x)) = B(x)F(\tfrac{1}{2}, \tfrac{2}{3}, \tfrac{1}{3}, 1, 1, C(x))\\
\text{3rd Class} \qquad & F(\tfrac{1}{2}, \tfrac{1}{2}, \tfrac{1}{2}, 1, 1, A(x)) = B(x)F(\tfrac{1}{2}, \tfrac{1}{4}, \tfrac{3}{4}, 1, 1, C(x))\\
\text{4th Class} \qquad & F(\tfrac{1}{2}, \tfrac{1}{2}, \tfrac{1}{2}, 1, 1, A(x)) = B(x)F(\tfrac{1}{2}, \tfrac{5}{6}, \tfrac{1}{6}, 1, 1, C(x))\\
\text{5th Class} \qquad & F(\tfrac{1}{2}, \tfrac{1}{3}, \tfrac{2}{3}, 1, 1, A(x)) = B(x)F(\tfrac{1}{2}, \tfrac{2}{3}, \tfrac{1}{3}, 1, 1, C(x))\\
\text{6th Class} \qquad & F(\tfrac{1}{2}, \tfrac{1}{3}, \tfrac{2}{3}, 1, 1, A(x)) = B(x)F(\tfrac{1}{2}, \tfrac{3}{4}, \tfrac{1}{4}, 1, 1, C(x))\\
\text{7th Class} \qquad & F(\tfrac{1}{2}, \tfrac{1}{3}, \tfrac{2}{3}, 1, 1, A(x)) = B(x)F(\tfrac{1}{2}, \tfrac{1}{6}, \tfrac{5}{6}, 1, 1, C(x))\\
\text{8th Class} \qquad & F(\tfrac{1}{2}, \tfrac{3}{4}, \tfrac{1}{4}, 1, 1, A(x)) = B(x)F(\tfrac{1}{2}, \tfrac{3}{4}, \tfrac{1}{4}, 1, 1, C(x))\\
\text{9th Class} \qquad & F(\tfrac{1}{2}, \tfrac{3}{4}, \tfrac{1}{4}, 1, 1, A(x)) = B(x)F(\tfrac{1}{2}, \tfrac{1}{6}, \tfrac{5}{6}, 1, 1, C(x))\\
\text{10th Class}\qquad  & F(\tfrac{1}{2}, \tfrac{1}{6}, \tfrac{5}{6}, 1, 1, A(x)) = B(x)F(\tfrac{1}{2}, \tfrac{1}{6}, \tfrac{5}{6}, 1, 1, C(x))
\end{align*}
where $A(x)$, $B(x)$, $C(x)$ are in general any functions and differ from class to class, of course.\\

And before we can use the translation method also practical, we will have to determine the transformations between the certain hypergeometric series. And below it will be seen more clarly, that this is the actual problem to solve. Because, the more transformations we have, the more often we can use the translation method. And this will also be the reason, which will force us to stop our investigations. But all this will become clearer.\\

To prove and/or find transformation formulas, it will be helpful, to reduce everything to the series $_2F_1$, because this one is easier to handle, because the number of variables is smaller. And for this we need a formula, relating the functions $_3F_2$ and $_2F_1$. And such a relation was given by T. Clausen in the memoir "Über die Fälle, wenn die Reihe $y = 1+\frac{\alpha}{1}\frac{\beta}{\gamma}x+\frac{\alpha\cdot(\alpha+1)}{1\cdot2}\cdot\frac{\beta(\beta+1}{\gamma(\gamma+1)}x^2+\text{etc.}$ ein Quadrat von der Form $z = 1+\frac{\alpha^{\prime}}{1}\cdot\frac{\beta^{\prime}}{\gamma^{\prime}}\cdot\frac{\delta^{\prime}}{\epsilon^{\prime}}x+\frac{\alpha^{\prime}(\alpha^{\prime}+1)}{1\cdot2}\cdot\frac{\beta^{\prime}(\beta^{\prime}+1)}{\gamma^{\prime}(\gamma^{\prime}+1)}\cdot\frac{\delta^{\prime}(\delta^{\prime}+1)}{\epsilon^{\prime}(\epsilon^{\prime}+1)}x^2+\text{etc.}$ hat" \cite{10}. There he states the following theorem
\subsection*{Theorem}
We have the identity
\[
	{}_2F_1^2(a, b, a+b+\tfrac{1}{2}, x) = {}_3F_2(2a, 2b, a+b, a+b+\tfrac{1}{2}, 2a+2b, x)
\]
This formula does not need to be proved here in all detail, it suffices - because our main focus is on the translation method itself - to explain the idea, how a proof could be given. The idea is the following one. At first the general functions $_2F_1(a,b,c,x)$ and (in the case of the Clausen identity) $_3F_2(a^{\prime},b^{\prime},c^{\prime},d^{\prime},e^{\prime},z)$ are consideres. From the definion as a power series one finds, that $_2F_1(a,b,c,x)$ satisfies this differential equation given by Euler in \cite{14}:
\[
(x^3-x^2)_2F_1^{\prime\prime}+[(a+b+1)x^2-cx]_2F_1^{\prime}+abx_2F_1=0
\]
And $_3F_2(a^{\prime},b^{\prime},c^{\prime},d^{\prime},e^{\prime},z)$ satisfies this one
\[
(x^3-x^2)_3F_2^{\prime\prime\prime}+[(3+a^{\prime}+b^{\prime}+d^{\prime})-(1+c^{\prime}+e^{\prime})]_3F_2^{\prime\prime}+
\]
\[
[(1+a^{\prime}+b^{\prime}+d^{\prime}+a^{\prime}b^{\prime}+a^{\prime}d^{\prime}+b^{\prime}d^{\prime})x-c^{\prime}e^{\prime}]_3F_2^{\prime}+a^{\prime}b^{\prime}d^{\prime}_3F_2=0
\]
So, to prove the Clausen identity a posteriori, one only has to check, that both sides satisfy the same differential equation with the same initial conditions, which reduces the proof to a mere calculation. We want to omit this calculation for that reason. Because the idea is more important for us here and we can prove all identities between hypergeometric series in the manner - a posteriori, of course .\\

After having seen the Clausen identity, we can focus our investigation on the function $_2F_1(a,b,c,x)$ and the derive from transformation formaulas found for $_2F_1(a,b,c,x)$ such for $_3F_2$.\\

And such transformations for $_2F_1(a,b,c,x)$ will then have this form:
\[
_2F_1(a,b,c,A(x))=B(x)_2F_1(a^{\prime},b^{\prime},c^{\prime},C(x))
\] 
where $A(x)$ in the most cases will just be the function $x$. Then we will only have to check, if both sides satisfy the differential equation given by Euler in \cite{14}.\\

Such transformations have been given by many mathematicians, some with most different methods. Therefore we just want to use them, because we only would have  to check them via the differential equation. At first we want to give two transformations for $_2F_1$, afterwards in the same manner some more for $_3F_2$, that allow a translation between the individual classes, that we constituted earlier. Not for all classes we will find such transformations (at least not such ones, that lead to formulas from our list), what will be seen by considering the following list of transformations.

\section*{Examples for transformations of the hypergeometric series}

For the transformation of the series $_2F_1(a,b,c,x)$ we want to note the following two.
\subsection*{Theorem 1}
We have
\[
_2F_1(a,b,c,x)=(1-x)^{c-a-b}_2F_1(c-a,c-b,c,x)
\]
This transformation was given at first by Euler in \cite{14}.
\subsection*{Theorem 2}
We have
\[
_2F_1(a,b,c,x)=(1-x)^{-b}_2F_1(c-a,b,c,\frac{x}{x-1})
\]
This transformation is sometimes also attributed to Euler, for example by E. Kummer in his first memoir on the hypergeometric series \cite{10}. But it seems, that it was at first written down explicitly by J. Pfaff in his textbook "Disquisitiones analyticae maxime ad calculem integralem et doctrinam serierum pertinentes" \cite{23}, which is, why it is often also called Pfaff's transformation. If you apply Pfaff's transformation two times in a row, you get Euler's, which gives more reason to think, that Euler was not aware of Pfaff's transformation, because he would surely have noticed this fact then.\\

The following transformations all refer to the function $_3F_2$ and can all be obtained from the corresponding ones for $_2F_1$ by using the Clausen identity. We list some examples and left some spaces empty on purpose - as mentioned above -, where it was impossible for the author of this paper, to find a useful one for the further investigation. The transformations are listed in the 10 classes from above.

\begin{flalign*}
	& \textsc{1st Class} ~~ (\text{from $s=\frac{1}{2}$ to $s=\frac{1}{2}$ })\\
	& \quad F\left(\tfrac{1}{2}, \tfrac{1}{2}, \tfrac{1}{2}, 1, 1, 4x(1-x)\right) = (1-x)^{-\frac{1}{2}}F\left(\tfrac{1}{2}, \tfrac{1}{2}, \tfrac{1}{2}, 1, 1, \tfrac{x^2}{4x-4}\right) & \\
	& \quad F\left(\tfrac{1}{2}, \tfrac{1}{2}, \tfrac{1}{2}, 1, 1, \tfrac{-4x}{(1-x)^2}\right) = (1-x)^{\frac{1}{2}}F\left(\tfrac{1}{2}, \tfrac{1}{2}, \tfrac{1}{2}, 1, 1, \tfrac{x^2}{4x-4}\right) & \\
&\text{Both can be derived from the transformations for $_2F_1$ given by Kummer. See \cite{22}}.\\	
	& \textsc{2nd Class} ~~ (\text{from $s=\frac{1}{2}$ to $s=\frac{1}{3}$ }) \\
	& \quad \text{none found} & \\
	& \textsc{3rd Class} ~~ (\text{from $s=\frac{1}{2}$ to $s=\frac{1}{4}$ }) \\
	& \quad F\left(\tfrac{1}{2}, \tfrac{1}{2}, \tfrac{1}{2}, 1, 1, x\right) = (1-x)^{\frac{1}{2}}F\left((\tfrac{1}{2}, \tfrac{1}{4}, \tfrac{3}{4}, 1, 1, \tfrac{-4x}{(1-x)^2}\right) & \\
	&\text{The formulas needed for this are found in Kummer's paper again \cite{22}}.\\
	& \textsc{4th Class} ~~ (\text{from $s=\frac{1}{2}$ to $s=\frac{1}{6}$ }) \\
	& \quad F\left( \tfrac{1}{2}, \tfrac{1}{2}, \tfrac{1}{2}, 1, 1, x \right) = \left( 1 - \tfrac{x}{4}\right)^{-\frac{1}{2}} F\left( \tfrac{1}{2}, \tfrac{1}{6}, \tfrac{5}{6}, 1, 1, \tfrac{27x^2}{(4-x)^3} \right) & \\
	& \quad  F\left(\tfrac{1}{2}, \tfrac{1}{2}, \tfrac{1}{2}, 1, 1, x\right) = (1-4x)^{-\frac{1}{2}}F\left( \tfrac{1}{2}, \tfrac{1}{6}, \tfrac{5}{6}, 1, 1, \tfrac{27x}{(4x-1)^3}\right) & \\
	&\text{For this class the formulas follow from Goursat's ones in \cite{17}}\\
	& \textsc{5th Class} ~~ (\text{from $s=\frac{1}{3}$ to $s=\frac{1}{3}$ }) \\
	& \quad F\left( \tfrac{2}{3}, \tfrac{1}{2}, \tfrac{1}{3}, 1, 1, \tfrac{256x^3}{9(3+x)^4}\right) = \tfrac{3+x}{3(1+3x)}F\left( \tfrac{1}{3}, \tfrac{2}{3}, \tfrac{1}{2}, 1, 1, \tfrac{256x}{9(1+3x)^4}\right) & \\
&\text{This formula was given in \cite{25} in this form by M.D. Rogers.}\\
	& \textsc{6th Class} ~~ (\text{from $s=\frac{1}{3}$ to $s=\frac{1}{4}$ }) \\ 
	& \quad \text{none found} & \\
	& \textsc{7th class} ~~ (\text{from $s=\frac{1}{3}$ to $s=\frac{1}{6}$ })\\
	& \quad F\left( \tfrac{2}{3}, \tfrac{1}{3}, \tfrac{1}{2}, 1, 1, 4x(1-x)\right) = \left(1 - \tfrac{8x}{9}\right)^{-\frac{1}{2}}F\left(\tfrac{1}{2}, \tfrac{1}{6}, \tfrac{5}{6}, 1, 1, \tfrac{64x^3(1-x)}{(9-8x)^3}\right) & \\
&\text{One consult Kummer \cite{22} or Goursat \cite{17} again.}\\
	& \textsc{8th Class} ~~ (\text{von $s=\frac{1}{4}$ zu $s=\frac{1}{4}$ }) \\
	& \quad F\left(\tfrac{1}{4}, \tfrac{3}{4}, \tfrac{1}{2}, 1, 1, \tfrac{-108x}{(1-16x)^3}\right) = \tfrac{1-16x}{1-4x}F\left(\tfrac{1}{2}, \tfrac{1}{4}, \tfrac{3}{4}, 1, 1, \tfrac{108x^2}{(1-4x)^3}\right) \\
&\text{This formula is even given like this in \cite{25} by Rogers}.\\
\end{flalign*}

\begin{flalign*}
& \textsc{9th Class} ~~ (\text{from $s=\frac{1}{4}$ to $s=\frac{1}{6}$ }) \\ 
	& \quad \text{none found} & \\
	& \textsc{10th Class} ~~ (\text{from $s=\frac{1}{6}$ to $s=\frac{1}{6}$ }) & \\ 
	& \quad \text{none found} &
\end{flalign*}
After having given some transformation formulas, the only thing missing is, that we prove a starting formula. For that we will use some known values of the Gau\ss ian hypergeoemtric series, which were also already given by Gau\ss{} \cite{16} himself.
\section*{Proof of a starting formula, so that the translation method can be used afterwards}

For this we only need two formulas, on the one hand the reflexion formula for $\Gamma({x})$, given by Euler in the memoir "Evolutio formulae integralis $\int x^{f-1}\diff{d}x\ln^{\frac{m}{n}}{x}$ integratione a valore $x = 0$ ad $x =1$ extensa" \cite{12}. This formula says
\[
\Gamma{(x)}\Gamma{(1-x)}=\frac{\pi}{\sin{(\pi x)}}
\]
And in addition we need the formula, which is often called Gau\ss 's second summation theorem, given by Gau\ss{} in \cite{16}
\[
	_2F_1( a,b,\tfrac{a}{2} + \tfrac{b}{2} + \tfrac{1}{2}, \tfrac{1}{2} ) = \frac{\Gamma (\tfrac{1}{2}) \Gamma (\tfrac{a+b+1}{2})}{\Gamma (\tfrac{a}{2} + \tfrac{1}{2} )\Gamma (\tfrac{b}{2} + \tfrac{1}{2} )}
\]
And by using these two formulas we can state the following 

\subsubsection*{Theorem}
We have
\[
	\sum_{n=0}^{\infty} \sum_{k=0}^{n}\frac{(s)_n (1-s)_{k} (s)_{n-k} (1-s)_{n-k}}{(1)_k (1)_{n-k}}\frac{n\cdot\left(\tfrac{1}{2}\right)^n}{n! (n-k)!} = \frac{2\sin{s\pi}}{\pi}
\]
\subsubsection*{Proof}
At first we see, that we have by the definition of  $_2F_1$ - we will for the sake of brevity write $F$ here:
\[
F(s, 1-s, 1, x)=\sum_{n=0}^{\infty}\frac{(s)_n(1-s)_n}{(1)_n}\frac{x^n}{n!}
\]
if one considers the square of this series, we find by the Cauchy multiplication formula
\[
	F^2(s, 1-s, 1, x) = \sum_{n=0}^{\infty} \sum_{k=0}^{n}\frac{(s)_k (1-s)_k}{(1)_k}\frac{(s)_{n-k}(1-s)_{n-k}}{(1)_{n-k}}\frac{x^n}{n! (n-k)!}
\]
Now just apply the operator $x\frac{\diff{d}}{\diff{d}x}$ to both sides of the equation, the right-hand side gives:
\[
	\sum_{n=0}^{\infty} \sum_{k=0}^{n}\frac{(s)_k (1-s)_k}{(1)_k}\frac{(s)_{n-k}(1-s)_{n-k}}{(1)_{n-k}}\frac{n x^n}{n!(n-k)!}
\]
and on the left-hand side we find by using the relation
\[
\frac{\diff{d}}{\diff{d}x}F(a,b,c,x)=\frac{ab}{c}F(a+1,b+1,c+1,x)
\]
which follows directly from the of $F$ definition as a series,
\[
x\frac{\diff{d}}{\diff{d}x}F^2(s,1-s,1,x)=2x\cdot s(1-s)F(s,1-s,1,x)F(s+1,2-s,2,x)
\]
Now just put $x=\frac{1}{2}$ and by using Gau\ss's second summation theorem for this particular value
\begin{align*}
		& s(1-s)\cdot F(s, 1-s, 1, \tfrac{1}{2})  F\left( s+1, 2-s, 2, \tfrac{1}{2}\right) \\
		=& s(1-s) \cdot \frac{\gam{\tfrac{1}{2}}\gam{1}}{\gam{\tfrac{s}{2} + \tfrac{1}{2}}\gam{1 - \tfrac{s}{2}}}\cdot\frac{\gam{\tfrac{1}{2}}\gam{2}}{\gam{1+\tfrac{s}{2}}\gam{\tfrac{3}{2} - \tfrac{s}{2}}}.
\end{align*}
If one uses the fundamental property of the Gamma function $\gam{x+1} = x\gam{x}$ and Euler's formula $\gam{x}\gam{1-x} = \frac{\pi}{\sin{x}}$, we find
\begin{align*}
	&4\cdot\frac{\gam{\tfrac{1}{2}}\gam{1}}{\gam{\tfrac{s}{2} + \tfrac{1}{2}}\gam{1 - \tfrac{s}{2}}}\cdot\frac{\gam{\tfrac{1}{2}}\gam{2}}{\gam{\tfrac{s}{2}}\gam{\tfrac{1}{2} - \tfrac{s}{2}}} = 4\cdot\pi\frac{\sin{\frac{s\pi}{2}}}{\pi}\cdot\frac{\sin{( \tfrac{s}{2} + \tfrac{1}{2})}\pi}{\pi} \\[1mm]
	=&\frac{4\sin{\frac{s\pi}{2}}\cdot\cos{\frac{s\pi}{2}}}{\pi} = \frac{2\sin{s\pi}}{\pi}.
\end{align*}
Q.E.D.\\

In the last two steps we used the addition theorem and duplication formula for $\sin{(x)}$.\\
This theorem immedeatly gives
\subsubsection*{Corollarv 1}
For $s = \frac{1}{2}, \frac{1}{3}, \frac{1}{4}, \frac{1}{6}$ we find the formulas proved by Z.W. Sun in the paper "Some new series for $\frac{1}{\pi}$ and related congruences" \cite{27}, $1.1 - 1.4$, namely:
\begin{flalign*}
	&\sum_{n=0}^{\infty} \frac{n}{32^n}\sum_{k=0}^{n}\binom{2k}{k}^2 \binom{2n-2k}{n-k}^2 = \frac{2}{\pi} && \\
	&\sum_{n=0}^{\infty} \frac{n}{54^n}\sum_{k=0}^{n}\binom{2k}{k}\binom{3k}{k}\binom{2(n-k)}{n-k}\binom{3(n-k)}{n-k} = \frac{\sqrt{3}}{\pi} && \\
	&\sum_{n=0}^{\infty} \frac{n}{128^n}\sum_{k=0}^{\infty}\binom{4k}{2k}\binom{2k}{k}\binom{4(n-k)}{2(n-k)}\binom{2(n-k)}{n-k} = \frac{\sqrt{2}}{\pi} && \\
	&\sum_{n=0}^{\infty} \frac{n}{864^n}\sum_{k=0}^{\infty}\binom{6k}{3k}\binom{3k}{k}\binom{6(n-k)}{3(n-k)}\binom{3(n-k)}{n-k} = \frac{1}{\pi}
\end{flalign*}

To see this, one just has to rewrite the occuring binomial coefficients. We have for example
\[
\frac{1}{32^n}\binom{2k}{k}^2 \binom{2n-2k}{n-k}^2=\frac{1}{2^n}\frac{(\frac{1}{2})_k^2(\frac{1}{2})_{n-k}^2}{(1)_k(1)_{n-k}}\frac{1}{n!(n-k)!}
\]
And similary we find the other three following formulas
\subsubsection*{Corollary 2}
Of course we could put in more values for $s$, to get more formulas for $\frac{1}{\pi}$.

So $s = \frac{1}{5}$ gives
\[
	\sum_{n=0}^{\infty} \sum_{k=0}^{n}\binom{-\frac{1}{5}}{k}\binom{-\frac{4}{5}}{k}\binom{-\frac{1}{5}}{n-k}\binom{-\frac{4}{5}}{n-k}n\left(\frac{1}{2}\right)^n = \sqrt{\frac{1}{2}(5 - \sqrt{5})}\cdot\frac{1}{\pi} = \frac{c}{\pi}
\]
But in this case $c^2$ is already not rational anymore, so that this formula does not satisfy the required conditions and is therefore not useful for the translation method, to which we will proceed, after all the preparations, now.
\section*{Application of the translation method}
At first we want to prove Bauer's series, as the following series, proved by G. Bauer in the treatise "Von der Koeffizienten der Reihen von Kugelfunktionen einer Variablen" \cite{3} is now called.
\subsubsection*{Theorem}
The following formula holds
\[
	\sum_{n=0}^{\infty} \frac{\pham{\frac{1}{2}}_n \pham{\frac{1}{2}}_n \pham{\frac{1}{2}}_n }{(1)_n^3}(4n+1)(-1)^n = \frac{2}{\pi}.
\]
\subsubsection*{Proof}
Consider Pfaff's transformation
\[
	_2F_1(\alpha, \beta, \gamma, x) = (1-x)^{-\beta} {}_2F_1(\gamma -\alpha, \beta, \gamma, \tfrac{x}{x-1})
\]
Apply this to the following square, where $\beta = \frac{1}{4}$. And as a consequence $\gamma= 1$ and $\alpha = \frac{3}{4}$.
\[
	F^2 (\tfrac{1}{4}, \tfrac{3}{4}, 1, x)
\]
so this gives
\[
_2F_1^2(\tfrac{1}{4},\tfrac{3}{4},1,x)=(1-x)^{-\tfrac{1}{2}} {}_2F_1^2(\tfrac{1}{4},\tfrac{1}{4},1,\tfrac{x}{x-1})
\]
Now you can use the Clausen identity on the right-hand side, which leads to
\[
_2F_1^2(\tfrac{1}{4},\tfrac{3}{4},1,x)=(1-x)^{-\tfrac{1}{2}} {}_3F_2(\tfrac{1}{2},\tfrac{1}{2},\tfrac{1}{2},1,1,\tfrac{x}{x-1})
\]
Apply the operator $x\frac{\diff{d}}{\diff{d}x}$ to both sides, the left-hand side is clear and becomes
\[
	\sum_{n=0}^{\infty} n\sum_{k=0}^{n}\frac{(\frac{1}{4})_k (\frac{3}{4})_k}{(1)_k}\frac{(\frac{1}{4})_{n-k}(\frac{3}{4})_{n-k}}{(1)_{n-k}}\frac{ x^n}{n!(n-k)!}
\]
the right-hand side, if one writes $a_n$ for the coefficient for the sake of brevity
\begin{align*}
	&x\frac{\diff{d}}{\diff{d}x}(1-x)^{-\tfrac{1}{2}} \sum_{n=0}^{\infty}a_n(\tfrac{x}{x-1})^n  \\
	=& x[\frac{1}{2}(1-x)^{-\frac{3}{2}}\sum_{n=0}^{\infty}a_n \pham{\tfrac{x}{x-1}}^n + (1-x)^{-\frac{1}{2}}\sum_{n=0}^{\infty}a_n \pham{\tfrac{x}{x-1}}^{n-1}(-\tfrac{1}{1-x)^2})] \\
	=&\frac{x}{2}(1-x)^{-\frac{3}{2}}\sum_{n=0}^{\infty}a_n \pham{\tfrac{x}{x-1}}^n+(1-x)^{-\frac{3}{2}} \sum_{n=0}^{\infty}na_n \pham{\tfrac{x}{x-1}}^n
\end{align*}
So we have this equation after all
\[
	\sum_{n=0}^{\infty} n\sum_{k=0}^{n}\frac{(\frac{1}{4})_k (\frac{3}{4})_k}{(1)_k}\frac{(\frac{1}{4})_{n-k}(\frac{3}{4})_{n-k}}{(1)_{n-k}}\frac{ x^n}{n!(n-k)!}
	= (1-x)^{-\frac{3}{2}}[\frac{x}{2}\sum_{n=0}^{\infty}a_n \pham{\tfrac{x}{x-1}}^n +\sum_{n=0}^{\infty}n a_n\pham{\tfrac{x}{x-1}}^2]
\]
For $x=\frac{1}{2}$ the left-hand side can be summed by the means of our theorem and we find - if we put $s=\frac{1}{4}$ there - this formula
\[
	\frac{\sin{\frac{\pi}{4}}}{\pi} = \frac{\sqrt{2}}{\pi} = \frac{\sqrt{2}}{2}\sum_{n=0}^{\infty}  \frac{\pham{\frac{1}{2}}_n^3}{(1)_n^3}(4n+1)(-1)^n
\]
and so
\[
	\frac{2}{\pi} = \sum_{n=0}^{\infty}  \frac{\pham{\frac{1}{2}}_n^3}{(1)_n^3}(4n+1)(-1)^n
\]
Q.E.D.\\

And this is already Bauer's series and the first example for the translation method. The Bauerian proof differs from the one given here, also the Ramanujanian one bases on different theorems, so that we have given a new proof for this series. In similar manner, just with another transformation for $_2F_1$, we find the following 
\subsubsection*{Theorem}
We have
\[
	\sum_{n=0}^{\infty}  \frac{\pham{\frac{1}{2}}_n^3}{(1)_n^3}(6n+1)(-\frac{1}{8})^n= \frac{2\sqrt{2}}{\pi}
\]
\subsubsection*{Proof}
We use Kummer's transformation from \cite{22}, which gives by squaring
\[
	_2F_1^2 (\tfrac{1}{2}, \tfrac{1}{2}, 1, x) = (1-x)^{-\frac{1}{2}}{} _2F_1^2(\tfrac{1}{4}, \tfrac{1}{4}, 1, \tfrac{x^2}{4x-4})
\]
On the right-hand side we can use Clausen's fomula again and our equation becomes
\[
	\sum_{n=0}^{\infty} \sum_{k=0}^{n}\frac{(\frac{1}{2})_k (\frac{1}{2})_k}{(1)_k}\frac{(\frac{1}{2})_{n-k}(\frac{1}{2})_{n-k}}{(1)_{n-k}}\frac{ x^n}{n!(n-k)!}=(1-x)^{-\frac{1}{2}} \sum_{n=0}^{\infty}  \frac{\pham{\frac{1}{2}}_n^3}{(1)_n^3}(\tfrac{x^2}{4x-4})^n
\]
For the sake of brevity we want to write $a_n$ for the coefficient on the left-hand side, on the right-hand side $b_n$ and for $\tfrac{x^2}{4x-4}$ $A(x)$, so that we will have
\[
\sum_{n=0}^{\infty} a_nx^n = (1-x)^{-\frac{1}{2}} \sum_{n=0}^{\infty} b_n(A(x))^n
\]
And now we proceed in the same way as just before, by applying the operator $x\frac{\diff{d}}{\diff{d}x}$ to both sides of the equations, so that we have
\[
\sum_{n=0}^{\infty} na_nx^n = \frac{x}{2}(1-x)^{-\frac{3}{2}}\sum_{n=0}^{\infty} b_n(A(x))^n -(1-x)^{-\frac{3}{2}}(x-2)\sum_{n=0}^{\infty} nb_n(A(x))^n
\] 
Now we set $x=\frac{1}{2}$ again, so that we can sum the left-hand side with our starting formula. This gives
\[
\frac{2}{\pi} = \frac{\sqrt{2}}{2}\sum_{n=0}^{\infty}  \frac{\pham{\frac{1}{2}}_n^3}{(1)_n^3}(-\tfrac{1}{8})^n+3\sqrt{2}\sum_{n=0}^{\infty}  \frac{\pham{\frac{1}{2}}_n^3}{(1)_n^3}n(-\tfrac{1}{8})^n
\]
and simplyfies to
\[
\frac{2\sqrt{2}}{\pi} = \sum_{n=0}^{\infty}  \frac{\pham{\frac{1}{2}}_n^3}{(1)_n^3}(6n+1)(-\tfrac{1}{8})^n
\]
Q.E.D.

And this is the second formula on our list.\\

It is already clear, that in the same way we can prove more formulas from our starting formula, as long as we do not run out of transformation formulas. Because the calculation is always the same, it suffices, to give the transformation for $_2F_1$. So we have the following
\subsection*{Theorem}
The following formula holds
\[
	\sum_{n=0}^{\infty}\frac{\pham{\frac{1}{6}}_n \pham{\frac{5}{6}}_n \pham{\frac{1}{2}}_n}{(1)_n^3} (28n + 3) \pham{\frac{27}{125}}^n = \frac{5\sqrt{5}}{\pi}
\]
\subsubsection*{Proof}
Just use this transformation given by Goursat
\[
_2F_1(\tfrac{1}{4}, \tfrac{3}{4}, 1, x) = (1-\tfrac{3x}{4})^{-\frac{1}{4}} {}_2F_1(\tfrac{1}{12}, \tfrac{5}{12}, 1, \tfrac{27x^2 (1-x)}{(4-3x)^3})
\]
this one just needs to be squared, so that we can use the Clausen identity on the right-hand side, which leads to
\[
_2F_1^2(\tfrac{1}{4}, \tfrac{3}{4}, 1, x)=(1-\tfrac{3x}{4})^{-\frac{1}{2}} {}_3F_2(\tfrac{1}{6}, \tfrac{5}{6},\tfrac{1}{2}, 1, 1, \tfrac{27x^2 (1-x)}{(4-3x)^3})
\] 
Then just apply the operator $x\frac{\diff{d}}{\diff{d}x}$ to both sides of the equation and put $x= \frac{1}{2}$ afterwards. Q.E.D. \\

In the same way we find
\subsubsection*{Theorem}
The following formula holds
\[
	\sum_{n=0}^{\infty} \frac{\pham{\frac{1}{6}}_n \pham{\frac{5}{6}}_n \pham{\frac{1}{2}}_n }{(1)_n^3}(22n + 2)\pham{\frac{4}{5^3}} = \frac{5\sqrt{15}}{3\pi}
\]
\subsubsection*{Proof}
Just proceed as in the proofs above and use this transformation of Goursat
\[
	_2F_1(\tfrac{1}{3}, \tfrac{2}{3}, 1, x) = (1-\tfrac{8x}{9})^{-\frac{1}{2}}{} _2F_1 (\tfrac{1}{12}, \tfrac{5}{12}, 1, \tfrac{64x^3 (1-x)}{(9-8x)^3})
\]
Q.E.D.

So we already obtained 4 series from our starting formula and proved them by elementary means.\\

We want to add another examples, which will serve as a warning, but also to explain a phenomenon.
\subsubsection*{Warning example}
Consider Gousat's transformation 
\[
	_2F_1(\tfrac{1}{3}, \tfrac{2}{3}, 1, x) = (1+8x)^{-\frac{1}{4}} {}_2F_1(\tfrac{1}{12}, \tfrac{5}{12}, 1, \tfrac{64x (1-x)^3}{(1+8x)^3}).
\]
Now we want to square both sides and use the Clausen identity on the right-hand side. This leads to
\[
	_2F_1^2 (\tfrac{1}{3}, \tfrac{2}{3}, 1, x) = (1+8x)^{-\frac{1}{2}} {}_3F_2(\tfrac{1}{2}, \tfrac{1}{6}, \tfrac{5}{6}, 1, 1, \tfrac{64x (1-x)^3}{(1+8x)^3}).
\]
If we consider the value $x = \frac{1}{2}$ we find
\[
	\frac{64x(1-x)^3}{(1+8x)^3} = \frac{4}{125}
\]
which value we find on our list and we already proved the corresponding series.If we introduce the abbreviations $a_n$, $b_n$ for the coefficients and the function $A(x)$ for the argument, this will, by applying the operator $x\frac{\diff{d}}{\diff{d}x}$ to both sides, on the right-hand side lead to
\begin{align*}
& x\frac{\diff{d}}{\diff{d}x} (1+8x)^{\frac{1}{2}}\sum_{n=0}^{\infty}b_n \pham{A(x)}^n \\
=& x\left\{ -\tfrac{8}{2}\cdot (1+8x)^{-\nfrac{3}{2}}\sum_{n=0}^{\infty} b_n (A(x))^n \right. \\
 &- \left. \frac{64(x-1)^2 (8x^2 + 20x - 1)}{(8x+1)^4}(1+8x)^{-\nfrac{1}{2}} \sum_{n=0}^{\infty} n b_n (A(x))^{n-1} \right\} \\
=& -4x(1+8x)^{-\nfrac{3}{2}}\sum_{n=0}^{\infty} b_n (A(x))^{n} \\
&- \frac{64x(x-1)^3}{(8x+1)^3}\frac{(8x^2 + 20x - 1)}{(8x+1)(x-1)}(1+8x)^{-\frac{1}{2}}\sum_{n=0}^{\infty}n b_n (A(x))^{n-1} \\
=& -4x(1+8x)^{-\frac{3}{2}}\sum_{n=0}^{\infty} b_n (A(x))^n + \frac{8x^2 + 20x - 1}{x-1}(1+8x)^{-\frac{3}{2}}\sum_{n=0}^{\infty} n b_n (A(x))^n
\end{align*}
So this is the right-hand side of our equation, the left-hand side can be summed by the means of our starting formula
\[
	\frac{\sqrt{3}}{\pi} = -\frac{2}{5\sqrt{5} }\sum_{n=0}^{\infty} b_n \pham{\tfrac{4}{125}}^n - \frac{22}{5\sqrt{5}}\sum_{n=0}^{\infty} n b_n \pham{\tfrac{4}{125}}^n
\]
or
\[
	-\frac{5\sqrt{5}}{\pi} = \sum_{n=0}^{\infty} (22n+2)b_n \pham{\tfrac{4}{125}}^n
\]
which is plainly wrong, because it differs from the previous formula.\\

So this example shows, that one always has to be careful and has to watch out for such things. The explaination of this phenomenon exceeds the elementary methods, that we used here.\\

We just want to note, that you can explain this result and this was for example done by A. Meurman in the appendix to the paper "Proof of some conjectured formulas for $\frac{1}{\pi}$ by Z.-W. Sun" \cite{1} .\\

But in the following we will avoid such cases, but it nevertheless appeared to be a good idea, to note this particular example, because it illustrates well, that those transformations have to be interpretated correctly. This will become more evident, when we get to the divergent series.\\

But now we want to prove some more formulas with the translation method, because we now can use a bit more formulas. They, for example, lead to the following

\subsubsection*{Theorem}
The two following summations are equivalent
\[
	\sum_{n=0}^{\infty}\frac{\pham{\frac{1}{2}}^3}{(1)_n^3}(6n+1)\pham{\frac{1}{4}}^n = \frac{4}{\pi}  
\]
and
\[
\sum_{n=0}^{\infty}\frac{\pham{\frac{1}{6}}_n \pham{\frac{5}{6}}_n \pham{\frac{1}{2}}_n}{(1)_n^3}(28n+3) = \frac{5\sqrt{5}}{\pi}.
\]
Here "equivalent" is to be understood, that the right equation follows from the left and vice versa, if the translation method is used.
\subsubsection*{Proof}
One just has to use the transformation
\[
	\pham{1 - \frac{x}{4}}^{-\frac{1}{2}} \sum_{n=0}^{\infty}\frac{\pham{\frac{1}{6}}_n \pham{\frac{5}{6}}_n \pham{\frac{1}{2}}_n}{(1)_n^3} \pham{\frac{27x^2}{(4-x)^3}}^n = \sum_{n=0}^{\infty}\frac{\pham{\frac{1}{2}}_n \pham{\frac{1}{2}}_n \pham{\frac{1}{2}}_n}{(1)_n^3}x^n.
\]
By an easy and analog calculation to the one above - so we will omit the calculation - it becomes clear, that by starting from the first equation you arrive at the second and vice versa by starting from the second at the first. Q.E.D.\\

Because we already proved the second equation, we also proved the second equation. And in this manner we can show more series. So we find
\subsubsection*{Theorem}
The following two formulas are equivalent:
\[
\sum_{n=0}^{\infty}\frac{\pham{\frac{1}{2}}_n \pham{\frac{1}{4}}_n \pham{\frac{3}{4}}_n}{(1)_n^3}(7n+1)\pham{\frac{32}{81}}^n = \frac{9}{2\pi} 
\]
and
\[
\sum_{n=0}^{\infty}\frac{\pham{\frac{1}{2}}_n^3}{(1)_n^3}(6n+1)\pham{-\frac{1}{8}}^n = \frac{2\sqrt{2}}{\pi}
\]
\subsubsection*{Proof}
Just translate one formula into the other by using the following translation
\[
	_3F_2(\tfrac{1}{2}, \tfrac{1}{2}, \tfrac{1}{2}, 1, 1, x) = (1-x)^{-\frac{1}{2}}{}_3F_2(\tfrac{1}{2}, \tfrac{1}{4}, \tfrac{3}{4}, 1, 1, -\tfrac{4x}{(1-x)^2} )
\]
Q.E.D.\\

And like this we find
\subsubsection*{Theorem}
The two following series are equivalent:
\[
\sum_{n=0}^{\infty}\frac{\pham{\frac{1}{2}}_n \pham{\frac{1}{6}}_n \pham{\frac{5}{6}}_n}{(1)_n^3}(256n + 20)\pham{\frac{2}{11}}^{3n} = \frac{11\sqrt{33}}{\pi}
\]
and
\[
\sum_{n=0}^{\infty}\frac{\pham{\frac{1}{2}}_n^3}{(1)_n^3}(6n+1)\pham{-\frac{1}{8}}^n = \frac{\sqrt{2}2}{\pi}
\]
\subsubsection*{Beweis}
Just translate one formula into the other by using this transformation
\[
	\sum_{n=0}^{\infty}\frac{\pham{\frac{1}{2}}_n^3 }{(1)_n^3}x^n = \pham{1-\frac{x}{4}}^{-\frac{1}{2}}\sum_{n=0}^{\infty}\frac{\pham{\frac{1}{2}}_n \pham{\frac{1}{6}}_n \pham{\frac{5}{6}}_n}{(1)_n^3}\pham{\frac{27x^2}{(4-x)^3}}^n
\]
Q.E.D.\\

So we got 7 series from our list and proved them completely. This might seem to be a small number, but one has to be aware of the fact, that the method, we explainend and used here, is more elementary than Ramanujan's and the WZ-Algorithm and therefore almost as a logical consequence can not lead that far as more sophisticated ones.\\

The idea of the translation method was now well illustrated on these examples and it becomes clear, that everything in the case of these formulas for $\frac{1}{\pi}$ reduces to finding the corresponding transformation formulas for the function $_3F_2$.\\

And this task, the finding of the transformations, gets more difficult as well. To get an impression of this, just consult Zudilins paper "Lost in Translation" \cite{31}, where the translation method is explained on one example. This example gives the following
\subsection*{Theorem}
The following two summations are equivalent:
\[
\sum_{n=0}^{\infty} \frac{(\frac{1}{2})_n^3 }{(1)_n^3}(4n + 1)\pham{-1}^n = \frac{2}{\pi}
\]
and
\[
\sum_{n=0}^{\infty} \frac{(\frac{1}{2})_n (\frac{3}{4})_n (\frac{1}{4})_n}{(1)_n^3}(40n + 3)\pham{\frac{1}{49^2}}^n = \frac{49\sqrt{3}}{9\pi}
\]
The corresponding transformation is really a monster and shows, how laborous the finding of the transformations becomes. And - as also shown in the mentioned paper \cite{31} - the work can be simplyfied, if one uses more sophisticated theories, in this case the theory of the modular functions.\\

So we want to assume now, that we already know Ramanujan's formulas, to see, how many formulas we can add.

\subsection*{Proving more formulas assuming Ramanujan's ones}
Becuase now we know the series
\[
	\sum_{n=0}^{\infty}\frac{\pham{\frac{1}{2}}_n^3}{(1)_n^3}(42n+5) \pham{\frac{1}{64}}^n = \frac{16}{\pi}
\]
we can prove another one from our list. But independent from this we would have the following theorem
\subsubsection*{Theorem}
The two following summations are equivalent
\[
\sum_{n=0}^{\infty}\frac{\pham{\frac{1}{2}}_n^3}{(1)_n^3}(42n+5) \pham{\frac{1}{64}}^n = \frac{16}{\pi}
\]
and
\[
\sum_{n=0}^{\infty}\frac{\pham{\frac{1}{2}}_n \pham{\frac{1}{4}}_n \pham{\frac{3}{4}}_n}{(1)_n^3}(65n+8)\pham{-\frac{16^2}{63^2}}^n
\]
\subsubsection*{Proof}
Use the following transformation for the translation
\[
	_3F_2(\tfrac{1}{2}, \tfrac{1}{2}, \tfrac{1}{2}, 1, 1, x) = (1-x)^{-\frac{1}{2}}{} _3F_2(\tfrac{1}{4}, \tfrac{3}{4}, \tfrac{1}{2}, 1, 1, -\tfrac{4x}{(1-x)^2})
\]
Q.E.D.

And with the same series we find
\subsubsection*{Theorem}
The following two summations are equivalent
\[
\sum_{n=0}^{\infty}\frac{\pham{\frac{1}{2}}_n^3}{(1)_n^3}(42n+5)\pham{\frac{1}{64}}^n = \frac{16}{\pi}
\]
and
\[
\sum_{n=0}^{\infty}\frac{\pham{\frac{1}{2}}_n \pham{\frac{1}{6}}_n \pham{\frac{5}{6}}_n}{(1)_n^3}(63n+8)\pham{-\frac{4}{5}}^{3n} = \frac{5\sqrt{15}}{\pi}
\]
\subsubsection*{Proof}
Just use the transformation
\[
	_3F_2(\tfrac{1}{2}, \tfrac{1}{2}, \tfrac{1}{2}, 1, 1, x) = (1 - 4x)^{-\frac{1}{4}} {}_3F_2(\tfrac{1}{2}, \tfrac{5}{12}, \tfrac{1}{12}, 1, 1, \tfrac{27x}{(4x-1)^3})
\]
Q.E.D.

And a third one can be proved

\subsubsection*{Theorem}
The two following series are equivalent
\[
	\sum_{n=0}^{\infty}\frac{\pham{\frac{1}{2}}_n^3}{(1)_n^3}(42n + 5)\frac{1}{64^n} = \frac{16}{\pi} \\
\]
and
\[
	\sum_{n=0}^{\infty}\frac{\pham{\frac{1}{2}}_n \pham{\frac{1}{6}}_n \pham{\frac{5}{6}}_n}{(1)_n^3}(2394n + 144)\pham{\frac{4}{85}}^{3n} = \frac{85\sqrt{85}}{\sqrt{3}\pi}
\]
\subsubsection*{Proof}
Just use this transformation formula
\[
	_3F_2(\tfrac{1}{2}, \tfrac{1}{2}, \tfrac{1}{2}, 1, 1, x) = (1-4x)^{-\frac{1}{2}} {}_3F_2(\tfrac{1}{2}, \tfrac{1}{6}, \tfrac{5}{6}, 1, 1, \tfrac{27x}{(4x-1)^3}).
\]
Q.E.D.

From this we see, that Ramanujan could have added these three series to his list - alogside with Bauer's series. Also in the book series, published by B. Berndt and G. Andrews "Ramanujan's Notebooks" \cite{2} you do not find this formulas, although there are many transformation formulas for the hypergeometric series - especially in the second one.\\

By using another series, proved by the means of the WZ-Algorithm, we can prove one more series, because we have the following
\subsection*{Theorem}
The two following series are equivalent:
\[
	\sum_{n=0}^{\infty}\frac{\pham{\frac{2}{3}}_n \pham{\frac{1}{3}}_n \pham{\frac{1}{2}}_n}{(1)_n^3}(5n+1)\pham{-\frac{9}{16}}^n = \frac{4}{\sqrt{3}\pi}
\]
and
\[
	\sum_{n=0}^{\infty}\frac{\pham{\frac{1}{6}}_n\pham{\frac{5}{6}}_n\pham{\frac{1}{2}}_n}{(1)_n^3}(506n + 31)\pham{-\frac{9}{40^3}}^n = \frac{160\sqrt{30}}{9\pi}.
\]
\subsubsection*{Beweis}
Because this proof is slightly different and to see at least one full calculation, we present it in all detail. We use the transformation
\[
	_3F_2(\tfrac{2}{3}, \tfrac{1}{3}, \tfrac{1}{2}, 1, 1, 4x(1-x)) = (1 - \tfrac{8x}{9})^{-\frac{1}{2}}{} _3F_2(\tfrac{1}{6}, \tfrac{5}{6}, \tfrac{1}{2}, 1, 1, \tfrac{64x^3 (1-x)}{(9-8x)^3}).
\]
We introduce $a_n$ and $b_n$ for the coefficients. Instead of $4x(1-x)$ we write $A(x)$, instead of $\frac{64x^3 (1-x)}{(9-8x)^3}$ $B(x)$. Then we will have
\[
	\sum_{n=0}^{\infty} a_n (A(x))^n = (1-\tfrac{8x}{9})^{-\frac{1}{2}}\sum_{n=0}^{\infty} b_n (B(x))^n.
\]
Now we apply the operator $4x(1-x)\frac{\diff{d}}{\diff{d}x}$ to both sides of the equation, on the left-hand side we find
\begin{align*}
& 4x(1-x)\frac{\diff{d}}{\diff{d}x}\sum_{n=0}^{\infty}a_n (A(x))^n \\
=& 4x(1-x)(4-8x)\sum_{n=0}^{\infty} n a_n (A(x))^{n-1} \\
=& (4-8x) \sum_{n=0}^{\infty} n a_n (A(x))^n
\end{align*}
And on the right-hand side
\begin{align*}
& 4x(1-x)\frac{\diff{d}}{\diff{d}x} \pham{1-\frac{8x}{9}}^{-\frac{1}{2}}\sum_{n=0}^{\infty} b_n (B(x))^n \\
=& 4x(1-x) \left\{ \frac{1}{2}\cdot\frac{8}{9}\pham{1 - \frac{8x}{9}}^{-\frac{3}{2}}\sum_{n=0}^{\infty} b_n (B(x))^n \right. \\
&+ \left. \pham{1-\frac{8}{9}}^{-\frac{1}{2}}\frac{64x^2(8x^2 - 36x + 27)}{(9-8x)^4}\sum_{n=0}^{\infty}n b_n (B(x))^{n-1} \right\} \\
=& \frac{16}{9}x(1-x)\pham{1 - \frac{8x}{9}}^{-\frac{3}{2}}\sum_{n=0}^{\infty} b_n (B(x))^n \\
&+ \pham{1 - \frac{8x}{9}}^{-\frac{1}{2}}\frac{48x^2 - 36x + 27}{9-8x}\sum_{n=0}^{\infty} n b_n (B(x))^n
\end{align*}
Now we consider the formula
\[
	\alpha\sum_{n=0}^{\infty} a_n (A(x))^n + \beta (4-8x)\sum_{n=0}^{\infty} n (A(x))^n
\]
and with it we want to obtain the formula
\[
	\sum_{n=0}^{\infty} a_n (5n+1) (-\tfrac{9}{16})^n.
\]
Therefore we solve $A(x) = 4x(1-x) = -\frac{9}{16}$ and $x = -\frac{1}{8}$ or $x = \frac{9}{8}$; we do not take the value $\frac{9}{8}$ here, because it would lead to a limit formula, which we will consider later. And the value $x=-\frac{1}{8}$ leads to
\[
	\alpha\sum_{n=0}^{\infty} a_n (A(x))^n + 5\beta\sum_{n=0}^{\infty} n(A(x))^n,
\]
so we have to choose $\alpha = 1$ and $\beta = 1$. And with the equations, we derived before, this leads to
\[
	\frac{160\sqrt{30}}{9\pi} = \sum_{n=0}^{\infty} (506n + 31)b_n (-\tfrac{9}{40^3})^n.
\]
Q.E.D.

So we got 4 more series and proved them, because we derived them form ones known from another source with translation.\\

We want to note, that Guillera independently by the means of the translation method proved new series from those, that he got with the WZ-method, how the author was told by him in a private communication, but he did not publish his results about this.\\

Now we want to proceed to divergent series of the same shape.
\section*{On divergent series of the same shape}
The following series have the same form as the previous ones and can be found on our list. And we will now see, how they arise from coverging ones via the translation method. We state the following 
\subsubsection*{Theorem}
We have
\[
	\sum_{n=0}^{\infty}\frac{\pham{\frac{1}{2}}_n}{(1)_n^3}(6n+1)(-8)^n = \frac{1}{\pi}.
\]
\subsubsection*{Proof}
We see, that this series diverges, but it can be proved with the translation method. Just consider Kummer's transformation
\[
_2F_1(\tfrac{1}{2}, \tfrac{1}{2}, 1, x) = (1-x)^{-\frac{1}{2}}{} _3F_2(\tfrac{1}{4}, \tfrac{1}{4}, 1, \tfrac{-4x}{(1-x)^2})
\]
by squaring and using the Clausen identity
\[
_2F_1^2(\tfrac{1}{2}, \tfrac{1}{2}, 1, x) = (1-x)^{-1}{}_3F_2(\tfrac{1}{2}, \tfrac{1}{2}, \tfrac{1}{2}, 1, 1, \tfrac{-4x}{(1-x)^2})
\]
Now we apply - as above - the operator $x\frac{\diff{d}}{\diff{d}x}$ to both sides and put $x=\frac{1}{2}$. The left-hand side can be summed by the means of our starting formula and we find calculating formally
\[
	\sum_{n=0}^{\infty} \frac{\pham{\frac{1}{2}}_n^3}{(1)_n^3}(6n+1)(-8)^n = \frac{1}{\pi}.
\]
Q.E.D.

And like this we find the next divergent series
\subsubsection*{Theorem}
We have
\[
	\sum_{n=0}^{\infty} \frac{\pham{\frac{1}{2}}_n \pham{\frac{3}{4}}_n \pham{\frac{1}{4}}_n}{(1)_n^3}(5n+1)(-\tfrac{16}{9})^n = \frac{\sqrt{3}}{\pi}.
\]
\subsubsection*{Proof}
Use the transformation
\[
	_3F_2(\tfrac{1}{2}, \tfrac{1}{2}, \tfrac{1}{2}, 1, 1, x) = (1-x)^{-\frac{1}{2}} {}_3F_2(\tfrac{1}{2}, \tfrac{1}{4}, \tfrac{3}{4}, 1, 1, -\tfrac{4x}{(1-x)^2}) 
\]
and the proved series
\[
\sum_{n=0}^{\infty} \frac{(\frac{1}{2})_n^3 }{(1)_n^3}(6n + 1)\pham{\frac{1}{4}}^n = \frac{4}{\pi}
\]
and translate it. Q.E.D.

By the same procedure we find
\subsubsection*{Theorem}
We habe formally
\[
	\sum_{n=0}^{\infty}\frac{\pham{\frac{1}{2}}_n \pham{\frac{1}{2}}_n \pham{\frac{1}{2}}_n}{(1)_n^3}(3n+1)4^n = -\frac{2\mathrm{i}}{\pi}.
\]
\subsubsection*{Proof}
This might be the weirdest formula up to now, because adding only positive terms leads to a complex number 

But we just do the translation as usual, we use this transformation
\[
	_3F_2(\tfrac{1}{2}, \tfrac{1}{2}, \tfrac{1}{2}, 1, 1, x) = (1-4x)^{-\frac{1}{2}}{} _3F_2(\tfrac{1}{2}, \tfrac{1}{6}, \tfrac{5}{6}, 1, 1, \tfrac{27x}{(4x-1)^3}).
\]
Now we let $\frac{27x}{(4x-1)^3} = \frac{4}{125}$ and find $x=4$ as a solution. For the translation we use the formaula
\[
	\sum_{n=0}^{\infty} \frac{\pham{\frac{1}{2}}_n \pham{\frac{5}{6}}_n \pham{\frac{1}{6}}_n}{(1)_n^3} (22n+2)(\tfrac{4}{125})^n = \frac{5\sqrt{15}}{3\pi}
\]
that we already proved. And for $x=4$ we find the formula in the theorem. Q.E.D.

It should be noted, that this series in the same way can be proved from this transformation
\[
_3F_2\left(\tfrac{1}{2}, \tfrac{1}{2}, \tfrac{1}{2}, 1, 1, x\right) = (1-x)^{-\frac{1}{2}}{}_3F_2\left((\tfrac{1}{2}, \tfrac{1}{4}, \tfrac{3}{4}, 1, 1, \tfrac{-4x}{(1-x)^2}\right)
\]
Because if you solve the equation $-\frac{4x}{(1-x)^2}=-\frac{16}{9}$  - we proved the corresponding divergent series - you will find the solutions $x=\frac{1}{4}$ - from the corresponding formula we proved the mentionend divergent series - and $x=4$. With this result, two divergent series can be translated into each other.

But we advance to the next series
\subsubsection*{Theorem}
Formally the following formula holds
\[
	\sum_{n=0}^{\infty}\frac{\pham{\frac{1}{2}}_n \pham{\frac{1}{3}}_n \pham{\frac{2}{3}}_n}{(1)_n^3}(15n+4)(-4)^n = \frac{3\sqrt{3}}{\pi}.
\]
\subsubsection*{proof}
Use Rogers' transformation
\[
	_3F_2(\tfrac{2}{3}, \tfrac{1}{3}, \tfrac{1}{2}, 1, 1, \tfrac{256x^3}{9(3+x)^4} ) = \tfrac{3+x}{3(1+3x)}{}_3F_2(\tfrac{2}{3}, \tfrac{1}{3}, \tfrac{1}{2}, 1, 1, \tfrac{256x}{9(1+3x)^4})
\]
and the series
\[
	\sum_{n=0}^{\infty}\frac{\pham{\frac{1}{2}}_n \pham{\frac{2}{3}}_n \pham{\frac{1}{3}}_n}{(1)_n^3}(4 + 33n)(\tfrac{4}{125})^n = \frac{15\sqrt{3}}{2\pi}.
\]
Q.E.D.

It will be convenient to mention, that Guillera used the divergent series, to prove the other one, how he shows in a talk given by him. He proved corresponding divergent series along with other divergent ones in his paper "WZ-proofs of "divergent" Ramanujan-type series via the WZ-method" \cite{19}. \\

And we further find
\subsubsection*{Theorem}
The following formula holds
\[
	\sum_{n=0}^{\infty}\frac{\pham{\frac{1}{2}}_n^3}{(1)_n^3}(21n+8)(64)^n = \frac{2\mathrm{i}}{\pi}.
\]
\subsubsection*{Proof}
For this one has to use the transformation
\[
	_3F_2(\tfrac{1}{2}, \tfrac{1}{2}, \tfrac{1}{2}, 1, 1, x) = (1-\tfrac{x}{4})^{-\frac{1}{2}}{} _3F_2(\tfrac{1}{2}, \tfrac{1}{6}, \tfrac{5}{6}, 1, 1, \tfrac{27x^2}{(4-x)^3})
\]
and the series
\[
	\sum_{n=0}^{\infty}\frac{\pham{\frac{1}{2}}_n \pham{\frac{1}{6}}_n \pham{\frac{5}{6}}_n}{(1)_n^3}(-\tfrac{4^3}{5^3})^n (63n+8) = \frac{5\sqrt{15}}{\pi}.
\]
Just note, that one solution of the equation $\frac{27x^2}{(4-x)^3}=-\frac{4^3}{5^3}$ is $x=64$. Q.E.D.

Of course, as also Guillera and Zudilin show in their joint paper "Divergent Ramanujan-type super congruences" \cite{21} on an example, divergent series can lead (and do often lead) to complex series as
\[
	\sum_{n=0}^{\infty}\frac{\pham{\frac{1}{2}}_n^3}{n!^3}\pham{\frac{105-21\sqrt{7}\mathrm{i}}{32}n + \frac{49 - 13\sqrt{7}\mathrm{i}}{64}}\pham{\frac{47 + 45\sqrt{7}\mathrm{i}}{128}}^n = \frac{\sqrt{7}}{\pi}.
\]
The authors say, that they got there, by considering the series
\[
	\sum_{n=0}^{\infty}\frac{\pham{\frac{1}{2}}_n^3}{(1)_n^3}(3n+1)2^{2n} = -\frac{2\mathrm{i}}{\pi}.
\]
This one can also be shown with our last theorem and by using this transformation
\[
	_3F_2(\tfrac{1}{2}, \tfrac{1}{2}, \tfrac{1}{2}, 1 ,1, 4x(1-x)) = (1-x)^{-\frac{1}{2}}{}_3F_2(\tfrac{1}{2}, \tfrac{1}{2}, \tfrac{1}{2}, 1, 1, \tfrac{x^2}{4x-4}).
\]
And more directly you get there with this transformation
\[
	_3F_2(\tfrac{1}{2}, \tfrac{1}{2}, \tfrac{1}{2}, 1, 1, x) = (1-x)^{-\frac{1}{2}}{} _3F_2(\tfrac{1}{2}, \tfrac{1}{6}, \tfrac{5}{6}, 1, 1, \tfrac{27x^2}{(4-x)^3})
\]
For this solve the equation $\frac{27x^2}{(4-x)^3} = -\frac{4^3}{5^3}$ and use the corresponding series from our list and finally note, that $\frac{47 + 45\sqrt{7}\mathrm{i}}{128}$ is a solution.

As in the paper by Guillera and Zudilin \cite{21} this gives the following
\subsubsection*{Theorem}
We have
\[
	\sum_{n=0}^{\infty}\frac{\pham{\frac{1}{2}}_n \pham{\frac{1}{4}}_n \pham{\frac{3}{4}}_n}{(1)_n^3}(30n+8)(\tfrac{4}{3})^{4n} = -\frac{18\mathrm{i}}{\pi}.
\]
\subsubsection*{Proof}
Translate the formula
\[
	\sum_{n=0}^{\infty}\frac{\pham{\frac{1}{2}}_n^3}{(1)n^3}\pham{\frac{105 - 21\sqrt{7}\mathrm{i}}{32}n + \frac{49 - 13\sqrt{7}\mathrm{i}}{64}}\pham{\frac{47 + 45\sqrt{7}\mathrm{i}}{128}} = \frac{\sqrt{7}}{\pi}
\]
with the transformation
\[
	{}_3F_2(\tfrac{1}{2}, \tfrac{1}{2}, \tfrac{1}{2}, 1, 1, z) = (1-x)^{-\frac{1}{2}} {}_3F_2(\tfrac{1}{2}, \tfrac{1}{4}, \tfrac{3}{4}, 1, 1, -\tfrac{4x}{(1-x)^2}).
\]
Q.E.D.

Now we saw some examples of divergent series of the same shape, which also have their meaning, if one takes into account, that we obtained them from the translation method. Therefore it can be interpretated, that you can get the convergent ones from the divergent ones.\\

And you can legitimate their use, if you just say, that those series have to be understood as analytic continuation of the function $_3F_2$.\\

But it seemed to be a good idea, to mention this after the formal proof of the formulas, to illustrate, that divergent series and formal calculation should not be abandoned already at the beginning. Euler by going through the divergent found remarkable results. In his textbook "Institutiones Calculi Differentilis" \cite{11} he gets to the fundamental properties of the $\Gamma-$ and $\Psi-$ (logarithmic derivate of the $\Gamma-$function) and discovers the functional equation of the Riemann zeta function in the paper "Remarques sur un beau rapport entre les series des puissances tant directes que reciproques" \cite{13}, to name just two things. And one rightly claims, that the theory of Fourier series would not be there, where it is today without divergent series and their formal use, if one takes into account, how often J. Fourier, the father of the theory of the doctrine of Fourier series, in his monumental book "Théorie analytique de la chaleur" \cite{15} gives proofs and does examples, involving divergent series. For example the first proof for the Fourier coefficients uses divergent series.\\

Why divergent series for $\frac{1}{\pi}$ often lead to complex values, can also be seen. The transformations often have this form
\[
_3F_2(a,b,\frac{1}{2},1,1,x)=\sqrt{1-(A(x))}{}_3F_2F(a^{\prime},b^{\prime},\frac{1}{2},1,1,B(x))
\]
where $A(x)$ in the most cases is just $x$ and $B(x)$ is a rational function. So if $A(x) > 1$, what happens rather often for divergent series, the transformation becomes imaginary. Just by accident, the imaginary quantities will cancel out. So divergent series have quite a potential to lead to complex series, as above.\\

We will not pursuit this any further, but want to mention, that such complex series were at first considered in the paper "Complex series for $\frac{1}{\pi}$" \cite{8} by Zudilin, H. Chan, J. Wan, but they obtain them by another method.\\

We now want to leave those examples, because we will not find more formulas from our list with the translation method - we just ran out $_3F_2$ of useful transformations for $_3F_2$ to do the translation. So the real problem consists - as mentionend earlier - in finding those transformations.\\

And after everything said and seen there will be likely noboby to doubt, that you from every series on the list by finding a certain transformation for the function $_3F_2$ a trasnlation into every other on the list can be done. This can be also formulated as follows.\\

All formulas on the list are equivalent by translation.\\

To this point this is only a conjeture, of which a proof is desired, because the way to the proof will without doubt bring along many nice results.\\

If this conjecture is true, you would only need one formula from this list to prove all the other in the same way we did it.\\

But for now we want to advance to formulas of the same kind, that arise by considering limits.

\section*{Limit formulas for $\frac{1}{\pi}$ arising from the translation method}
The first theorem will provide an example, that will show the nature of these formulas
\subsubsection*{Theorem}
The following limits hold
\[  \lim_{x \rightarrow \frac{1}{2}} \frac{1-2x}{1-x}\sum_{n=0}^{\infty} n\frac{\pham{\frac{1}{2}}_n \pham{\frac{1}{2}}_n \pham{\frac{1}{2}}_n}{(1)_n^3}(4x(1-x))^n = \frac{2}{\pi} \]
\[  \lim_{x \rightarrow \frac{1}{2}} \frac{1-2x}{1-x}\sum_{n=0}^{\infty} n \frac{\pham{\frac{1}{2}}_n \pham{\frac{1}{3}}_n \pham{\frac{1}{2}}_n}{(1)_n^3}(4x(1-x))^n = \frac{\sqrt{3}}{\pi} \] 
\[ \lim_{x \rightarrow \frac{1}{2}} \frac{1-2x}{1-x}\sum_{n=0}^{\infty} n\frac{\pham{\frac{1}{2}}_n \pham{\frac{1}{4}}_n \pham{\frac{3}{4}}_n}{(1)_n^3}(4x(1-x))^n = \frac{\sqrt{2}}{\pi} \]
\[  \lim_{x \rightarrow \frac{1}{2}} \frac{1-2x}{1-x}\sum_{n=0}^{\infty} n\frac{\pham{\frac{1}{2}}_n \pham{\frac{5}{6}}_n \pham{\frac{1}{6}}_n}{(1)_n^3}(4x(1-x))^n = \frac{1}{\pi} \]
\subsubsection*{Proof}
Consider the square
\[
	_2F_1^2 (s, 1-s, 1, x)
\]
This can be transformed with Kummer's transformation
\[
	_2F_1 (s, 1-s, 1, x) = {}_2F_1 (\tfrac{s}{2}, \tfrac{1-s}{2}, 1, 4x(1-x))
\]
into
\[
_2F_1^2 (s, 1-s, 1, x) = {}_2F_1^2 (\tfrac{s}{2}, \tfrac{1-s}{2}, 1, 4x(1-x))
\]
Now you can use Clausen's identity
\[
	_2F_1^2 (s, 1-s, 1, x) = {}_3F_2(s, 1-s, \tfrac{1}{2}, 1, 1, 4x(1-x)). 
\]
Just apply the operator $x\frac{\diff{d}}{\diff{d}x}$ to both sides of the equation, the right-hand side with $a_n$ for the coefficient and $A(x)$ for $4x(1-x)$ becomes
\begin{align*}
	& x\frac{\diff{d}}{\diff{d}x}\sum_{n=0}^{\infty}\frac{(s)_n (1-s)_n (\frac{1}{2})_n}{(1)_n^3}(4x(1-x)) \\
	=& x\frac{\diff{d}}{\diff{d}x} \left\{ \sum_{n=0}^{\infty} a_n (A(x))^n \right\} \\
	=& x(4-8x)\sum_{n=0}^{\infty} n a_n (A(x))^{n-1} \\
	=& 4\frac{1-x}{1-x}x(1-2x)\sum_{n=0}^{\infty} n a_n (A(x))^{n-1} \\
	=& \frac{1-2x}{1-x} \sum_{n=0}^{\infty} n a_n (A(x))^n
\end{align*}
The left-hand side for $x=\frac{1}{2}$ can be summed with our starting formula, on the richt-hand side you have to take the limit for $x \rightarrow \frac{1}{2}$, this leads to the formula
\[
	\frac{2\sin{s\pi}}{\pi} = \lim_{x \rightarrow \frac{1}{2}} \frac{1-2x}{1-x}\sum_{n=0}^{\infty} n a_n (A(x))^n
\]
And for $s = \frac{1}{2}, \frac{1}{3}, \frac{1}{4}, \frac{1}{6}$ we obtain the values from the theorem. Q.E.D.\\

It is easyly seen, why one has to take limits, because the series
\[
\sum_{n=0}^{\infty} n a_n (A(x))^n
\]
diverges for $\frac{1}{2}$ and this is compensated by the value $1-2x$, that converges to zero. So the procedure - as it will also be in the following limit formulas - leads to limits of the kind "$0\cdot \infty$". We directly go on to the following theorems

\subsubsection*{Theorem}
We have the following limit
\[
	\lim_{x \rightarrow -\frac{1}{8}} (8x+1)\sum_{n=0}^{\infty}\frac{\pham{\nfrac{1}{2}}_n \pham{\nfrac{1}{6}}_n \pham{\nfrac{5}{6}}_n}{n!^3}n \pham{\tfrac{27x}{(4x-1)^3}}^n = \frac{\sqrt{3}}{2\pi}.
\]
\subsubsection*{Beweis}
Here we also use the translation method

\[
	\sum_{n=0}^{\infty}\frac{\pham{\frac{1}{2}}_n}{(1)_n^3} x^n = (1-4x)^{-\frac{1}{2}}\sum_{n=0}^{\infty}\frac{\pham{\frac{1}{2}}_n \pham{\frac{1}{6}}_n \pham{\frac{5}{6}}_n}{(1)_n^3}\pham{\frac{27x}{(4x-1)^3}}^n
\]
and the series
\[
	\sum_{n=0}^{\infty}\frac{\pham{\frac{1}{2}}_n^3}{(1)_n}(6n+1)\pham{-\tfrac{1}{8}}^n = \frac{\sqrt{2}\cdot 2}{\pi}.
\]
Because we have to be careful here, we want to display the important steps. As always we apply the operator $x\frac{\diff{d}}{\diff{d}x}$ to both sides, the left-hand side is clear, the right-hand side on the other hand leads to
\begin{align*}
	& x\frac{\diff{d}}{\diff{d}x} (1-4x)^{-\frac{1}{2}}\sum_{n=0}^{\infty}\frac{\pham{\frac{1}{2}}_n \pham{\frac{5}{6}}_n \pham{\frac{1}{6}}_n }{(1)_n^3}\pham{\tfrac{27x}{(4x-1)^3}}^n \\
	=& x\frac{\diff{d}}{\diff{d}x}(1-4x)^{-\frac{1}{2}}\sum_{n=0}^{\infty} a_n (A(x))^n \\
	=& x \left\{\tfrac{4}{2}(1-4x)^{-\frac{3}{2}}\sum_{n=0}^{\infty} a_n (A(x))^n - \tfrac{27(8x+1)}{(4x-1)^4}(1-4x)^{-\frac{1}{2}}\sum_{n=0}^{\infty}n a_n (A(x))^{n-1}\right\} \\
	=& 2x(1-4x)^{-\frac{3}{2}}\sum_{n=0}^{\infty} a_n (A(x))^n + \tfrac{27x}{(4x-1)^3}\tfrac{(8x+1)(1-4x)^{-\frac{1}{2}}}{(1-4x)^1}\sum_{n=0}^{\infty} n a_n (A(x))^{n-1} \\
	=& 2x(1-4x)^{-\frac{3}{2}}\sum_{n=0}^{\infty}a_n (A(x))^n + (8x+1)(1-4x)^{-\frac{3}{2}}\sum_{n=0}^{\infty}n a_n (A(x))^n
\end{align*}
and from this we derive the following general equation for
\[
	\sum_{n=0}^{\infty}\frac{\pham{\frac{1}{2}}_n^3}{(1)_n}(A+Bn)\pham{x}^n
\]
where $A$ and $B$ are arbitrary numbers
\begin{align*}
	=& A(1-4x)^{-\frac{1}{2}} \sum_{n=0}^{\infty}a_n (A(x))^n \\ 
	&+ B\left\{ 2x(1-4x)^{-\frac{3}{2}}\sum_{n=0}^{\infty}a_n (A(x))^n + (8x+1)(1-4x)^{-\frac{3}{2}}\sum_{n=0}^{\infty} n a_n (A(x))^n\right\} \\
	=& \left\{ A(1-4x)^{-\frac{1}{2}} + B\cdot 2x(1-4x)^{-\frac{3}{2}}\right\}\sum_{n=0}^{\infty} a_n (A(x))^n \\
	&+ B(8x+1)(1-4x)^{-\frac{3}{2}}\sum_{n=0}^{\infty}n a_n (A(x))^n.
\end{align*}
And for $x=-\frac{1}{8}$, $A=1$, $B=6$ we get a known series on the left-hand side and can use its value, on the right-hand side we have to take limits again. This leads to
\begin{align*}
	\frac{2\sqrt{2}}{\pi} =& \lim_{x\rightarrow -\frac{1}{8}}[\left\{1\cdot(1-4x)^{-\frac{1}{2}} + 6\cdot 2x(1-4x)^{-\frac{3}{2}}\right\}\sum_{n=0}^{\infty} a_n (A(x))^n \\
	&+ 6\cdot\tfrac{\sqrt{2}}{2}(8x+1)(1-4x)^{-\frac{3}{2}}\sum_{n=0}^{\infty}n a_n (A(x))^n].
\end{align*}
The first part of the limit involves the series $_3F_2(\tfrac{1}{6}, \tfrac{5}{6}, \tfrac{1}{2}, 1, 1; 1)$, which is easyly be seen to converge and this is everything we need here. And we find the limit
\[
	\frac{2\sqrt{2}}{\pi} = \lim_{x\rightarrow -\frac{1}{8}} 4\sqrt{\tfrac{2}{3}}(8x+1)\sum_{n=0}^{\infty}n a_n (A(x))^n
\]
or
\[
	\frac{\sqrt{3}}{2\pi} = \lim_{x\rightarrow -\frac{1}{8}}(8x+1)\sum_{n=0}^{\infty}n a_n (A(x))^n
\]
This proves the theorem. Q.E.D.

And like this you find the following
\subsubsection*{Theorem}
We have the following limit
\[
	\frac{\sqrt{2}}{\pi} = \lim_{x\rightarrow 1^-}(x+1)\sum_{n=0}^{\infty}n\frac{\pham{\frac{1}{2}}_n \pham{\frac{3}{4}}_n \pham{\frac{1}{4}}_n}{(1)_n^3}\pham{-\tfrac{4x}{(1-x)^2}}^n.
\]
\subsubsection*{Proof}
By using the transformation
\[
	_3F_2(\tfrac{1}{2}, \tfrac{1}{2}, \tfrac{1}{2}, 1, x) = (1-x)^{-\frac{1}{2}}{}_3F_2(\tfrac{1}{2}, \tfrac{1}{4}, \tfrac{3}{4}, 1, 1, -\tfrac{4x}{(1-x)^2})
\]
translate the series
\[
	\sum_{n=0}^{\infty}\frac{\pham{\frac{1}{2}}_n^3}{(1)_n^3}(4n+1)(-1)^n = \frac{2}{\pi}.
\]
Q.E.D.

In the same way you find
\subsubsection*{Theorem}
We have the limit
\[
	\frac{4\sqrt{3}}{\pi} = \lim_{x\rightarrow -8}\sum_{n=0}^{\infty}n\frac{\pham{\frac{1}{6}}_n \pham{\frac{5}{6}}_n \pham{\frac{1}{2}}_n}{(1)_n^3}\pham{\tfrac{27x^3}{(4-x)^3}}^n.
\]
\subsubsection*{Proof}
To prove this limit, translate the series
\[
	_3F_2(\tfrac{1}{2}, \tfrac{1}{2}, \tfrac{1}{2}, 1, 1, x) = (1-\tfrac{x}{4})^{-\frac{1}{2}}{}_3F_2(\tfrac{1}{6}, \tfrac{5}{6}, \tfrac{1}{2}, 1, 1, \tfrac{27x^2}{(4-x)^3})
\]
with the divergent series
\[
	\sum_{n=0}^{\infty}\frac{\pham{\frac{1}{2}}_n^3}{(1)_n^3}(3n+1)(-8)^n = \frac{1}{\pi}.
\]
The calculation is exactly like the one above and the limit is correct despite the divergent series. Q.E.D.

And like this it is possible, to find more limit formulas, like, for example from the tranformation, given by Rogers in \cite{25}.
\[
	_3F_2(\tfrac{1}{4}, \tfrac{1}{2}, \tfrac{3}{4}, 1, 1, \tfrac{256u^3}{9(3+u)^4}) = \tfrac{3+u}{3(1+3u)}{} _3F_2(\tfrac{1}{2}, \tfrac{1}{4}, \tfrac{1}{4}, 1, 1, \tfrac{256u}{9(1+3u)^4})
\]
and the formula
\[
	\sum_{n=0}^{\infty}\frac{\pham{\frac{1}{2}}_n \pham{\frac{1}{4}}_n \pham{\frac{3}{4}}_n}{(1)_n^3}(40n+3)\pham{\frac{1}{49^2}}^n = \frac{49\sqrt{3}}{9\pi}
\] 
But we want to leave it at this, because the subject of these limit formulas has hardly considered. Indeed such formulas seem to have their origin in a unpublished manuscript by the author of this paper last year - this unpublished manuscript is the one cited in Guillera's paper "Kinds o Proofs of Ramanujan-like series" \cite{20} - and here, where we only want to explain the translation method, is not the place, to explain these limit formulas in more detail. Nevertheless those formulas are worthy, to be be considered more extensively.\\

Because if they arise from the translation method, they can also be translated, to prove other series.\\

With this insight it might become easier to solve the formulated problem, to prove by the means of the translation method, that all formulas on our list are equivalent.\\

But now we will proceed to Suns's conjectures.

\section*{On Zhi-Wei Sun's conjectures}
In his paper "List of conjectured series for powers of $\pi$ and other constants" \cite{29} Sun gives many conjectures, which are divided into sections. The most conjectures in the fourth section were proved in the paper "Proof of some conjectured formulas for $\frac{1}{\pi}$ by Z.W. Sun" \cite{1}. Indeed only $4.14$ remains unproved, which we want to show now.
\subsection*{Theorem}
The following formula holds
\[
	\sum_{n=0}^{\infty}\frac{3n-1}{2^n}\sum_{k=0}^{n}\binom{-\frac{1}{3}}{k}\binom{-\frac{2}{3}}{n-k}\binom{-\frac{1}{6}}{k}\binom{-\frac{5}{6}}{n-k} = \frac{3\sqrt{6}}{2\pi}.
\]
\subsection*{Proof}
At first we see, that we have, by resolving the binomial coefficients
\[
	_2F_1(\tfrac{2}{6}, \tfrac{1}{6}, 1, x) = \sum_{n=0}^{\infty}\binom{-\frac{1}{3}}{n}\binom{-\frac{1}{6}}{n}x^n
\]
and
\[
	_2F_1(\tfrac{4}{6}, \tfrac{5}{6}, 1, x) = \sum_{n=0}^{\infty}\binom{-\frac{4}{6}}{n}\binom{-\frac{5}{6}}{n}x^n.
\]
Therefore the Cauchy product of the two series will give
\[
	_2F_1(\tfrac{1}{6}, \tfrac{2}{6}, 1, x){}_2F_1(\tfrac{4}{6}, \tfrac{5}{6}, 1, x) = \sum_{n=0}^{\infty}x^n\sum_{k=0}^{n}\binom{-\frac{1}{6}}{k}\binom{-\frac{2}{6}}{k}\binom{-\frac{4}{6}}{n-k}\binom{-\frac{5}{6}}{n-k}
\]
We want to transform the series
\[
_2F_1(\tfrac{4}{6}, \tfrac{5}{6}, 1, x)
\]
with Euler's transformation
\[
	F(\alpha, \beta, \gamma, x) = (1-x)^{\gamma - \alpha - \beta} F(\gamma - \alpha, \gamma - \beta, \gamma, x) 
\]
Then we will have for the product
\[
	_2F_1(\tfrac{1}{6}, \tfrac{2}{6}, 1, x){}_2F_1(\tfrac{4}{6}, \tfrac{5}{6}, 1, x) = (1-x)^{-\frac{1}{2}}{}_2F_1^2(\tfrac{1}{6}, \tfrac{2}{6}, 1, x).
\]
To the square the Clausen identity can be applied. This will lead to
\[
_2F_1(\tfrac{1}{6}, \tfrac{2}{6}, 1, x){}_2F_1(\tfrac{4}{6}, \tfrac{5}{6}, 1, x) = (1-x)^{-\frac{1}{2}}{}_3F_2(\tfrac{1}{3}, \tfrac{2}{3}, \tfrac{1}{2}, 1, 1, x)
\]
or
\begin{align*}
&\sum_{n=0}^{\infty}x^n\sum_{k=0}^{\infty}\binom{-\frac{1}{2}}{k}\binom{-\frac{2}{6}}{k}\binom{-\frac{4}{6}}{n-k}\binom{-\frac{5}{6}}{n-k} &=& (1-x)^{-\frac{1}{2}}\sum_{n=0}^{\infty}\frac{\pham{\frac{1}{2}}_n \pham{\frac{2}{3}}_n \pham{\frac{1}{3}}_n}{(1)_n^3}x^n 	\\
&\sum_{n=0}^{\infty}x^n c_n \hfill &=& (1-x)^{-\frac{1}{2}}\sum_{n=0}^{\infty} a_n x^n
\end{align*}
Now we can apply the operator $x\frac{\diff{d}}{\diff{d}x}$ to both sides again and proceed as usual. The left-hand side is clear
\[
	x\frac{\diff{d}}{\diff{d}x}\sum_{n=0}^{\infty}c_nx^n = \sum_{n=0}^{\infty}n c_n x^n.
\]
the right-hand side on the other hand is
\begin{align*}
&x\frac{\diff{d}}{\diff{d}x} (1-x)^{-\frac{1}{2}}\sum_{n=0}^{\infty} a_n x^n \\
=& x\left\{ \tfrac{1}{2}(1-x)^{-\frac{3}{2}} \sum_{n=0}^{\infty}a_n x^n + (1-x)^{-\frac{1}{2}} \sum_{n=0}^{\infty}n a_n x^{n-1}\right\} \\
=& \tfrac{x}{2}(1-x)^{-\frac{3}{2}}\sum_{n=0}^{\infty}a_n x^n + (1-x)^{-\frac{1}{2}}\sum_{n=0}^{\infty}n a_n x^n
\end{align*}
From this, with arbitrary $A$ and $B$ we deduce the general equation
\begin{align*}
\sum_{n=0}^{\infty} (A+Bn)c_n x^n &= A(1-x)^{-\frac{1}{2}}\sum_{n=0}^{\infty} a_n x^n \\
&+ B\left\{\tfrac{x}{2}(1-x)^{-\frac{3}{2}}\sum_{n=0}^{\infty}a_n + (1-x)^{-\frac{1}{2}}\sum_{n=0}^{\infty}n a_n x^n \right\} \\
&=\left\{ A(1-x)^{-\frac{1}{2}} + \tfrac{Bx}{2}(1-x)^{-\frac{3}{2}}\right\}\sum_{n=0}^{\infty} a_n x^n \\
&+ B(1-x)^{-\frac{1}{2}}\sum_{n=0}^{\infty} n a_n x^n
\end{align*}
Now we put $A=1$, $B=3$ and $x=\frac{1}{2}$ and find
\begin{align*}
\sum_{n=0}^{\infty}(3n+1)c_n\pham{\tfrac{1}{2}}^n 
=&\frac{\sqrt{2}}{2}\left\{ 6\sum_{n=0}^{\infty}(6n+1)a_n\pham{\tfrac{1}{2}}^n\right\}
\end{align*}
The formula on the right-hand side can be found on our list
\[
	\sum_{n=0}^{\infty}(6n+1)a_n\pham{\tfrac{1}{2}}^n = \frac{3\sqrt{3}}{\pi}
\]
and completes the proof. Q.E.D.

And in the same way we can add more formulas of this type to Sun's list.\\

Because we see, that all formulas from $4.2$ to $4.13$ have this form
\[
	\sum_{n=0}^{\infty}(\alpha n + \beta )x^n\sum_{k=0}^{n}\binom{2k}{k}\binom{-s}{k}\binom{2(n-k)}{n-k}\binom{-(1-s)}{n-k}
\]
with $s = \frac{1}{3}$ or $s = \frac{1}{4}$.

But because we also have Ramanujan-like formulas, that lead to $s = \frac{1}{2}$ and $s=\frac{1}{6}$ it will be worth the effort, to consider these ones.

We prove the following
\subsubsection*{Theorem}
We have the following general transformation
\[
	\sum_{n=0}^{\infty}x^n\binom{2k}{k}\binom{-s}{k}\binom{2(n-k)}{n-k}\binom{-(1-s)}{n-k} = {}_3F_2(\tfrac{1}{2}, s, 1-s, 1, 1, -\tfrac{4x^2}{1+4x}).
\]
\subsubsection*{Proof}
At first, we see, that the following formula holds
\[
	_2F_1(\tfrac{1}{2}, s, 1, -4x) = \sum_{k=0}^{\infty}\binom{2k}{k}\binom{-s}{k}x^k
\]
which follows from elementary identities for binomial coeficients, and if we write $1-s$ for $s$ 
\[
	_2F_1(\tfrac{1}{2}, 1-s, 1, -4x) = \sum_{n=0}^{\infty}\binom{2n}{n}\binom{-(1-s)}{n}x^n
\]
And by multiplying these two series
\[
	\sum_{n=0}^{\infty}x^n\sum_{k=0}^{n}\binom{2k}{k}\binom{-s}{k}\binom{2(n-k)}{n-k}\binom{-(1-s)}{n-k} = {}_2F_1(\tfrac{1}{2}, s, 1, -4x){}_2F_1(\tfrac{1}{2}, 1-s, 1, -4x)
\]
so we can consider hypergeometric series and their transformations again. Here we use Euler's transformation
\[
	_2F_1(\alpha, \beta, \gamma, x) = (1-x)^{\gamma - \alpha - \beta}{}_2F_1(\gamma - \alpha, \gamma - \beta, \gamma, x)
\]
which we want to apply to the second series, the product will then become
\begin{align*}
&_2F_1(\tfrac{1}{2}, s, 1, -4x){}_2F_1(\tfrac{1}{2}, 1-s, 1, 4x) \\
=&_2F_1(\tfrac{1}{2}, s, 1, -4x)(1+4x)^{1-\frac{1}{2}-(1-s)}{}_2F_1(\tfrac{1}{2}, s, 1, -4x) \\
=& (1+4x)^{-\frac{1}{2}+s}{}_2F_1^2(\tfrac{1}{2}, s, 1, -4x)
\end{align*}
Now we use Kummer's transformation
\[
	_2F_1(s, \tfrac{1}{2}, 1, -4x) = (1+4x)^{-\frac{s}{2}}{}_2F_1(\tfrac{s}{2}, \tfrac{1-s}{2}, 1, -\tfrac{4x^2}{1+4x})
\]
So we can continue the calculation and find
\begin{align*}
&(1+4x)^{-\frac{1}{2}+s}{}_2F_1^2(\tfrac{1}{2}, s, 1, -4x)\\
=& (1+4x)^{-\frac{1}{2}}{}_2F_1^2(\tfrac{s}{2}, \tfrac{1-s}{2}, 1, -\tfrac{4x^2}{1+4x})
\end{align*}
And with the Clausen identity
\[
	=(1+4x)^{-\frac{1}{2}}{}_3F_2(\tfrac{1}{2}, s, 1-s, 1,1, -\tfrac{4x^2}{1+4x})
\]
Q.E.D.

And from this we can deduce new formulas of the Sunian type, by using this theorem and the formulas from our list for $s=\frac{1}{2}$ and $s=\frac{1}{6}$ and translate them. We consider the case $s =\frac{1}{2}$ at first.
\subsection*{Corollary 1} 
We have
\[
	\sum_{n=0}^{\infty}n\pham{-\tfrac{1}{8}}^n\sum_{k=0}^{n}\binom{2k}{k}\binom{-\frac{1}{2}}{k}\binom{2(n-k)}{n-k}\binom{-\frac{1}{2}}{n-k} = \frac{4}{\sqrt{2}\pi}
\]
and
\[
	\sum_{n=0}^{\infty}(2n+1)\pham{\tfrac{1}{4}}^n\sum_{k=0}^{n}\binom{2k}{k}\binom{-\frac{1}{2}}{k}\binom{2(n-k)}{n-k}\binom{-\frac{1}{2}}{n-k} = \frac{4}{2\sqrt{2}\pi}.
\]
The proof is clear, just note, that $-\frac{4x^2}{1+4x}=-\frac{1}{8}$ has the two solutions $x = -\frac{1}{8}$ and $x = \frac{1}{4}$. Therefore one has just to translate the following series
\[
\sum_{n=0}^{\infty} \frac{(\frac{1}{2})_n^3 }{(1)_n^3}(6n + 1)\pham{-\frac{1}{8}}^n
\]
with our theorem.\\

We remark that the second formula appears in Suns's paper that "Some new series for $\frac{1}{\pi}$ and related congruences" \cite{27} as congruence. See $4.1$. There you find
\[
	\sum_{k=0}^{p-1}\frac{2n+1}{(-16)^n}\sum_{k=0}^{k}\binom{2k}{k}^2\binom{2n-2k}{n-k}^2\equiv p\pham{-\tfrac{1}{p}} + 3p^3 E_{p-3} \mod{p^4}
\]
- where one just has to rewrite the binomial coefficients. 
Here $p$ is an odd prime, $(\tfrac{p}{q})$ is the Legendre-Symbol and $E_n$ is the nth Euler number. 

And also the first formula of this corollary is just $1.1$ represented slightly different.\\

The connection between such formulas $\frac{1}{\pi}$ and such congruences is highly interesting, which is why we chose this example.

We want to consider the value $s=\frac{1}{6}$ now.

\subsubsection*{Corollary 2}
We have
\[
	\sqrt{3}\sum_{n=0}^{\infty} \pham{\frac{51}{250} + n\frac{192}{125}}\pham{\frac{3}{500}}^n\sum_{k=0}^{\infty}\binom{2k}{k}\binom{-\frac{1}{6}}{k}\binom{2n-2k}{n-k}\binom{-\tfrac{5}{6}}{n-k} = \frac{160\sqrt{30}}{9\pi}
\]
and
\[
	\sqrt{3}\sum_{n=0}^{\infty}\pham{\frac{21}{128} + n\frac{375}{256}}\pham{\frac{-3}{512}}^n\sum_{k=0}^{n}\binom{2k}{k}\binom{-\frac{1}{6}}{k}\binom{2n-2k}{n-k}\binom{-\frac{5}{6}}{n-k} = \frac{160\sqrt{30}}{9\pi}
\]
These summations follow from the series
\[
	\sum_{n=0}^{\infty}\frac{\pham{\frac{1}{6}}_n \pham{\frac{5}{6}}_n \pham{\frac{1}{2}}_n}{(1)_n^3}(506n + 31)\pham{-\frac{9}{40^3}}^n = \frac{160\sqrt{30}}{9\pi}
\]
by translation. The arguments explain themselves from the solutions of the equation $-\frac{4x^2}{1+4x}=\frac{9}{40^3}$ - $x = -\frac{3}{512}$ and $x = \frac{3}{500}$.\\

The formulas are presented here as in the joint paper with Almkvist "Proof of some conjectured formulas for $\frac{1}{\pi}$ by Zhi-Wei Sun" \cite{1}. We repeated the proof of the general theorem here, because it shows, that the translation method also gives some new formulas.\\

Of course here it helped a lot, to have the formulas and conjectures, whose consideration led to the general theorems. Therefore Sun's paper "Conjectures and results on $x^2 \mod{p^2}$ with $4p = x^2 + dy^2$" \cite{28} will helpful for everyone, who wants to discover new formulas for $\frac{1}{\pi}$.\\ 

And it becomes clear, that we could expand this Sunian list, because there are more transformations, that  provide the desired things. So we deduce this transformation from known transformations for $_2F_1$
\[
	_2F_1^2(\tfrac{1}{2}, \tfrac{1}{2}, 1, -4x) = (1+4x)^1 {}_3F_2(\tfrac{1}{2}, \tfrac{1}{2}, \tfrac{1}{2}, 1, 1, \tfrac{16x}{(1+4x)^2})
\]
which theorem is similar to that one, that we just gave..\\

Just note, that the left-hand side leads to the Sunian formulas, we considered here, for $s=\frac{1}{2}$. If we solve the equation $\frac{16z}{(1+4z)^2} = -8$ - we find the cooresponding divergent series on our list -, we find $x = -\frac{1}{2}$ and $x = -\frac{1}{8}$. This for $x = -\frac{1}{8}$ gives a already known formula, $x = -\frac{1}{2}$ would lead to a divergent series of this type. So there are some more formulas of this type to be found, because we just considered one other transformation for $_2F_1$, that was useful for this case.

But now we want to prove another one of Sun's conjectures 
\subsection*{Proof of Suns's formula 2.11}
We want to show
\[
\sum_{n=0}^{\infty}\frac{4k+1}{(-192)^k}\binom{2k}{k}S_k^{(2)}(4) = \frac{\sqrt{3}}{\pi}
\]
which is formula $2.11$ in Sun's Paper "List of conjectural series for powers of $\pi$ and other constants" \cite{29} - with
\[
	S^{(2)}_n (x) = \sum_{k=0}^{n}\binom{2k}{k}\binom{2n-2k}{n-k}x^{n-k}.
\]
For this we note at first, that we have
\subsection*{Theorem}
\[
	\binom{2n}{n}S_n^{(2)}(4) = \sum_{k=0}^{n}\binom{2k}{k}\binom{4k}{2k}\binom{4(n-k)}{2(n-k)}\binom{2(n-k)}{n-k}
\]
\subsection*{Proof}
Although noone will doubt this by comparing the formulas in Suns's paper "Some new series for $\frac{1}{\pi}$ and related congruences" \cite{27} - $(1.3)$ and $(1.11)$ - which are respectively
\[
	\sum_{n=0}^{\infty}\frac{n}{128^n}\sum_{k=0}^{n}\binom{2k}{k}\binom{4k}{2k}\binom{4(n-k)}{2(n-k)}\binom{2(n-k)}{n-k} = \frac{\sqrt{2}}{\pi},
\]
and
\[
	\sum_{n=0}^{\infty}\frac{n}{128^n}\binom{2n}{n}S_n^{(2)}(4) = \frac{\sqrt{2}}{\pi}
\]
we want to prove this. For this we note
\[
	\binom{2n}{n}S_n^{(2)}(4) = 4^n\binom{2n}{n}^2 {}_3F_2(\tfrac{1}{2}, \tfrac{1}{2}, -n, 1, \tfrac{1}{2}-n, 1)
\]
and
\[
	\sum_{k=0}^{n}\binom{2k}{k}\binom{4k}{2k}\binom{2(n-k)}{n-k}\binom{4(n-k)}{2(n-k)} = \binom{2n}{n}\binom{4n}{2n}{}_4F_3(\tfrac{1}{4}, \tfrac{3}{4}, -n, n, 1, \tfrac{1}{4}-n, \tfrac{3}{4}-n, 1)
\]
In L. Slater's book "Generalized Hypergeometric Functions"  \cite{26} we find this identity (Theorem 2.5.10)
\begin{align*}
	&{}_3 F_2(2a, -1, 2b, -n, 2c-1, \tfrac{1}{2}+a+b, -c-n, 1) \\
	=& \frac{(c-a+\tfrac{1}{2})_n (c-b-\tfrac{1}{2})_n}{(c-\tfrac{1}{2})_n (c+\tfrac{1}{2}-a-b)_n}{}_4 F_3 (a,b, 1-c-n, -n, \tfrac{1}{2}+a-c-n, \tfrac{3}{2} + b - c - n, 1)
\end{align*}
from which the theorem follows by subtituting the corresponding values. Q.E.D.\\

You can prove this formula from Sun's formulas (2.3) and (2.5) in \cite{27}, but Sun got these formulas with the WZ-method, that we want to avoid in this paper.

So this theorem allows us, to rewrite Suns's conjectured formula $2.11$ and we have

\subsection*{Corollary}
We have
\begin{align*}
\sum_{n=0}^{\infty} x^n S_n^{(2)}(4)\binom{2n}{n} &= \sum_{n=0}^{\infty}x^n\sum_{k=0}^{n}\binom{2k}{k}\binom{4k}{2k}\binom{2(n-k)}{n-k}\binom{4(n-k)}{2(n-k)} \\
&= {}_2F_1^2(\tfrac{1}{4}, \tfrac{3}{4}, 1, 64x)
\end{align*}
In addition we note, that by using Pfaff's transformation and the Clausen identity afterwards
\begin{align*}
	{}_2F_1^2(\tfrac{1}{4}, \tfrac{3}{4}, 1, x) =& ((1-x)^{-\frac{1}{4}})^2 {}_2F_1^2(\tfrac{1}{4},\tfrac{1}{4}, 1, \tfrac{x}{x-1}) \\
		=& (1-x)^{-\frac{1}{2}}{}_3F_2(\tfrac{1}{2}, \tfrac{1}{2}, \tfrac{1}{2}, 1, 1, \tfrac{x}{x-1})
\end{align*}
So we have 
\[
\sum_{n=0}^{\infty} x^n S_n^{(2)}(4)\binom{2n}{n}={}_3F_2(\tfrac{1}{2}, \tfrac{1}{2}, \tfrac{1}{2}, 1, 1, \tfrac{64x}{64x-1})
\]
and can prove the following theorem
\subsubsection*{Theorem}
Sun's conjectured formula
\[
	\sum_{k=0}^{\infty}\frac{4k+1}{(-192)^k}\binom{2k}{k}S_k^{(2)}(4) = \frac{\sqrt{3}}{\pi}
\]
is correct.
\subsubsection*{Proof}
Just translate the following formula with our theorem
\[
	\sum_{n=0}^{\infty}\frac{\pham{\frac{1}{2}}_n^3}{(1)_n^3}(6n+1)\pham{\tfrac{1}{4}}^n = \frac{4}{\pi}
\]
to give a proof. Q.E.D.

Of course we can deduce more formulas of this kind from our general theorem, if we use the corresponding formulas from our list .\\

Just solve the following equations and translate the corresponding formulas
\[
	\frac{64x}{64x-1} = -1,\quad -\frac{1}{8},\quad \frac{1}{64},\quad 4,\quad -8,\quad 64
\]
which respectively lead to formulas with the following values
\[
	\frac{1}{128},\quad -\frac{1}{576},\quad -\frac{1}{4032},\quad \frac{1}{48},\quad \frac{1}{72},\quad \frac{1}{63}
\]
The first three values were proved by Sun in \cite {27}, the value $\frac{1}{72}$ was covered by A. Meurman in \cite{1}. The last two values seem to be new.

Because the procedure is clear now, we want to conclude and give some final remarks.
\section*{Conclusion and final remarks}
So now we have explained the nature of the translation method in quite detail on the example of Ramanujan's formulas. And we only used rather elementary theorems, techniques and ideas.

And in addition it was shown, that some formulas for $\frac{1}{\pi}$ of this kind - and probably all - on our list are equivalent by translation.

And the translation method also leads to other kinds of formulas for $\frac{1}{\pi}$. Just consider the transformation given by Rogers in his paper "New ${}_5F_4$ hypergeometric transformations, three-variable Mahler measures, and formulas for $\frac{1}{\pi}$" \cite{25}
\[
	{}_3 F_2(\tfrac{1}{4}, \tfrac{1}{2}, \tfrac{3}{4}, 1, 1, \tfrac{256u}{9(1+3u)^4}) = \frac{1+3u}{1+u}\sum_{n=0}^{\infty}\pham{\tfrac{u}{9(1+u)^2}}^n\binom{2n}{n}\sum_{k=0}^{n}\binom{2k}{k}\binom{n}{k}^2
\]
Now we can use this formula from our list
\[
	\sum_{n=0}^{\infty}(40n+3)\pham{\frac{1}{2}}_n \pham{\frac{3}{4}}_n \pham{\frac{1}{4}}_n\pham{\frac{1}{49^2}}^n = \frac{2\sqrt{3}}{\pi}
\]
and translate it, to show
\[
	\sum_{n=0}^{\infty}(16n+3)\pham{\tfrac{1}{100}}^n\binom{2n}{n}\sum_{k=0}^{n}\binom{2k}{k}\binom{n}{k}^2 = \frac{25}{\sqrt{3}\pi}
\]
This formula at first was communicated to the author of this paper by Almkvist and there had its origins in the consideration of polytopes and K3 surfaces.\\

Zudilin, Chan, Y. Tanigawa und Y. Yang in the paper "New analogues of Clausen's identities arising from the theory of modular forms" \cite{6} - as the title suggests - deduce this one and many others of this kind from the theory of modular forms. We on the other hand found one example with translation.\\

Also formulas for $\frac{1}{\pi}$ of the kind
\[
	\sum_{n=0}^{\infty}\frac{(s)_n(1-s)_n}{(1)_n^2}P_n(x)z^n
\]
where $P_n(x)$ is the nth Legendre polynomial, have already been considered. In Sun's conjectures in \cite{29} these formulas occupy the whole fifth section - although there the connection to the Legendre polynomials is not immedeatly clear. This was at first shown in the paper "Legendre Polynomials and Ramanujan-type series for $\frac{1}{\pi}$" \cite{7} by Chan, Wan and Zudilin.\\

And considering all these papers and everything said, it is natural to conjecture, that all these kinds of formulas for $\frac{1}{\pi}$ - we just considered and metionend very few - can be translated into each other, which means that probably all of these formulas are equivalent to each other, although we still have no clear proof of this, which will likely require a theory, that explains this connections.\\
\section*{Acknowledgement}
The author would like to thank D. van Straten, G. Almkvist and W. Zudilin for their encouragement to write this paper and answering all the occuring questions. Also the author would like to thank J. Guillera and Z.-W. for the nice communication on the subject and making useful suggestions.\\

And last but not least the author would like to thank A. Diener, who helped very much with the \LaTeX{} version of this paper and completely selfless wrote the \LaTeX{} version of the preversion of this paper.

\end{document}